\pdfoutput=1
\documentclass[11pt,14paper]{amsart}
\usepackage[margin=17mm]{geometry}
\usepackage{amsmath,amsthm}
\usepackage{amssymb}

\usepackage[colorlinks]{hyperref}
\usepackage[numbers]{natbib}
\usepackage[all]{xypic}

\usepackage{graphicx}
\usepackage{epsfig}

\def\aabd#1 {\advance\abovedisplayskip by#1pt
\advance\belowdisplayskip by#1pt }
\def\abs#1 {\global\advance\baselineskip by#1pt}

\newtheorem{theorem}{Theorem}
\newtheorem{corollary}{Corollary}
\newtheorem{definition}{Definition}
\newtheorem{lemma}{Lemma}
\newtheorem{proposition}{Proposition}
\newcommand{\x}{\boldsymbol{x}}
\newcommand{\N}{\mathbb{N}}
\newcommand{\Z}{\mathbb{Z}}
\newcommand{\R}{\mathbb{R}}
\newcommand{\C}{\mathbb{C}}

\newcommand{\p}{\mathbb{P}}
\newcommand{\gh}{\mathfrak h}

\newcommand{\gp}{\mathfrak p}

\newcommand{\gl}{\mathfrak{gl}}

\newcommand{\Span}[1]{\left\langle#1\right\rangle}
\DeclareMathOperator{\LL}{LGr}

\DeclareMathOperator{\tr}{tr}

\DeclareMathOperator{\id}{id}
\DeclareMathOperator{\rank}{rank}
\DeclareMathOperator{\vol}{vol}

\DeclareMathOperator{\OOO}{\mathsf{O}}
\DeclareMathOperator{\CO}{\mathsf{CO}}

\DeclareMathOperator{\GL}{\mathsf{GL}}

\DeclareMathOperator{\Sp}{\mathsf{Sp}}

\DeclareMathOperator{\Gr}{Gr}

\DeclareMathOperator{\Fl}{Fl}%

\DeclareMathOperator{\Cont}{Cont}

\DeclareMathOperator{\Hom}{Hom}
\DeclareMathOperator{\hess}{hess}

\DeclareMathOperator{\Id}{Id}

\DeclareMathOperator{\Smbl}{Smbl}
\DeclareMathOperator{\Sing}{Sing}
\DeclareMathOperator{\dd}{d}
\DeclareMathOperator{\Stab}{\mathrm{Stab}}
\newcommand{\OO}{\mathcal{O}}
\newcommand{\E}{\mathcal{E}}
\newcommand{\FF}{\mathcal{F}}

\newcommand{\CC}{\mathcal{C}}
\newcommand{\cU}{\mathcal{U}}
\newcommand{\sR}{\mathcal{R}}
\newcommand{\ins}{\Bigg\lrcorner}
\def\sp{\mathfrak{sp}}

\def\iso{\mathrm{iso}}

\let\CC\CC

 \def\polhk#1{\setbox0=\hbox{#1}{\ooalign{\hidewidth
  \lower1.5ex\hbox{`}\hidewidth\crcr\unhbox0}}}

\renewcommand{\S}{Section\ }

\begin{document}



\title[Lagrangian Grassmannians and PDEs]{Geometry of Lagrangian Grassmannians and nonlinear PDEs}


\author{Jan Gutt}
\address{Center for Theoretical Physics of the Polish Academy of Sciences,\\
Al. Lotnik\'ow 32/46, 
02-668 Warsaw,
Poland\\
E-mail: jan.gutt@gmail.com}
\author{Gianni Manno}
\address{Dipartimento di Science Matematiche ``G. L. Lagrange'', Politecnico di Torino,\\Corso Duca degli Abruzzi, 24, 10129 Torino, Italy.\\
 E-mail:  giovanni.manno@polito.it}
 \author{Giovanni Moreno}
 \address{Department of Mathematical Methods in Physics,\\
 Faculty of Physics, University of Warsaw,\\
ul. Pasteura 5, 02-093 Warszawa, Poland\\
E-mail: giovanni.moreno@fuw.edu.pl}
\thanks{Research founded by the Polish National Science Centre grant under the contract number 2016/22/M/ST1/00542.}

\maketitle

\begin{abstract}
This paper contains a thorough introduction to the basic geometric properties of the manifold of Lagrangian subspaces of a linear symplectic space,
known as the Lagrangian Grassmannian. It also reviews the important relationship between hypersurfaces in the Lagrangian Grassmannian and second-order PDEs.
Equipped with a comprehensive bibliography, this paper has been especially designed as an opening contribution for the proceedings volume of the homonymous
workshop held in Warsaw, September 5--9, 2016, and organised by the authors.
\end{abstract}

\goodbreak
{\small
\setcounter{tocdepth}{2}
\tableofcontents

}
\abs-0.15 \aabd-1
\section*{Introduction}\addcontentsline{toc}{section}{Introduction}
One way to see how geometry enters  the theory of second order PDEs in one dependent and $n$ independent variables is to regard the Hessian matrix
$(p_{ij})$ as an $n$-dimensional subspace $\Span{e_i+p_{ij}\epsilon^j\mid i=1,\dots, n}$ in the space $\R^n\oplus\R^{n\ast}$. Obviously this correspondence
is not accidental: its motivations will be thoroughly surveyed in Section~\ref{secPart2}. However, a crucial fact can already be noticed: the symmetry of
the Hessian matrix corresponds to the subspace $\Span{e_i+p_{ij}\epsilon^j }$ being \emph{isotropic} with respect to the canonical symplectic form
on $\R^n\oplus\R^{n\ast}$. In other words, the aforementioned subspace is \emph{Lagrangian}, that is an element of the \emph{Lagrangian Grassmannian}
$\LL(n,2n)$.

The object $\LL(n,2n)$ exists independently of theories of  PDEs. It is indeed a very well-known projective variety, displaying a lot of interesting
properties, smoothness and homogeneity above all else. As such, it can be studied \emph{per se}, and this is precisely the purpose of
Section~\ref{secPart1}. Due to the algebro-geometric origin of $\LL(n,2n)$, we shall examine the real case (relevant for applications to PDEs) as much as
possible in parallel with the complex case.

The entire content of this paper can be found elsewhere scattered throughout the existing literature. Our main goal was to squeeze a plethora of tiny
small elements---facts, formulas, lemmas, properties, remarks---into a short self-contained introductory paper. Taken individually they may seem
trivial, but their appropriate combination against the motivating background of PDEs form an unexpectedly rich and coherent picture.

The present paper serves yet another purpose. It is included in the Banach Center Publications volume dedicated to the workshop titled
\textit{Geometry of Lagrangian Grassmannians and Nonlinear PDEs} and held in Warsaw in September 2016. The volume is designed in such a way as to
provide a source book for a monographic graduate/postgraduate course, as well as a reference for recent research in the discipline
(Section~\ref{secPart3}). The present paper may represent a good departing point for the novice. It may also guide the expert finding his/her way in
the rest of the volume (see Section~\ref{secPart3.1}).

\subsection*{Background, motivations and acknowledgements}
One of the main driving forces behind the present paper, the volume it belongs to and the homonymous workshop has been a conjecture, formulated
in 2010 by Ferapontov and his collaborators about the class of second order hydrodynamically integrable PDEs. Essentially, the conjecture states that
multidimensional hydrodynamically integrable second order PDEs  of Hirota type are of Monge--Amp\`ere type, see \cite[Section 1]{DoubrovFerapontov}.
Intrigued by Ferapontov's problem, two of us (GM and GM) started a systematic study of the notion of hydrodynamic integrability and soon realised
that there was a lot of differential and algebraic geometry at play. More complementary competences were needed. A first informal meeting was held
in 2014 in Milan, bringing the problem to the attention of Musso and Russo (both contributors to this volume). Interesting links with the geometry
of special projective varieties and homogeneous spaces were highlighted. In 2015 one of us (Moreno) was granted a two-year Maria Sk\l odowska-Curie
Fellowship at IMPAN (Warsaw) for continuing the study of the geometry of hypersurfaces in the Lagrangian Grassmannian and second order PDEs. It was
during this period that the authors of the present paper began their cooperation. In 2016 they organised the aforementioned workshop and started
editing the present volume.

To date, the conjecture is still open, even though it triggered an enormous amount of side and related works, eventually leading to interesting independent
results. The authors wish first of all to thank Professor Ferapontov for his deep and insightful analysis of the phenomenon of hydrodynamic integrability
and regret he could not make it to a workshop built, in a sense, around a his idea. Many thanks go also to all the other speakers and contributors
to this volume, to Professors: Hwang for his surprise visit, Bryant and Ciliberto for important remarks and valuable advices. The authors thank also
Professor Mormul  for reviewing the manuscript.

The authors acknowledge  the support of the Maria Sk\l odowska-Curie fellowship SEP-210182301 ``GEOGRAL'', the Institute of Mathematics of the Polish
Academy of Sciences, the Banach Centre, the project
 ``FIR (Futuro in Ricerca) 2013 -- Geometria delle equazioni differenziali'',   the grant 346300 for IMPAN from the Simons Foundation and the matching
2015--2019 Polish MNiSW fund. Giovanni Moreno has been also partially founded by the  Polish National Science Centre grant
under the contract number 2016/22/M/ST1/00542. Gianni Manno was partially supported by a ``Starting Grant per Giovani Ricercatori'' 53\_RSG16MANGIO
of the Polytechnic of Turin. Gianni Manno and Giovanni Moreno are members of G.N.S.A.G.A of I.N.d.A.M.\newpage

  \section{Geometry of the (real and complex) Lagrangian Grassmannian}\label{secPart1}
\subsection{Preliminaries}\label{secPrelim}
One of the harshest lessons from earlier studies in Mathematics is the impossibility to identify a vector space with its dual in a canonical way. This is
mirrored in Physics by the profound difference between vectors and covectors. The former correspond geometrically to those fancy arrows emanating from 0,
whereas   the latter are \emph{hyperplanes} passing through 0.

However, if a ``balanced mixture'' of vectors and covectors is given, such as in the space $V\oplus V^*$, then there is an obvious way to identify the space
with its dual. Just perform a ``counterclockwise rotation by $\frac{\pi}{2}$'', having identified the horizontal axis with~$V$ and the vertical axis with
$V^*$. The evident analogy with the multiplication by $i$ in the complex plane led to the coinage of the term \emph{symplectic} by Hermann Weyl in 1939
\cite[page 165]{Weyl1939}. Indeed the preposition ``sym'' is the Greek analog of the Latin preposition ``cum'', see, e.g., \cite[pp.~xiii--xiv]{daSilva2008}
and \cite{MO45162}.

From now on, $V$ is a linear vector space of dimension $n$, and $\omega$ is the 2-form on $V\oplus V^*$ corresponding to the  canonical identification
$(V\oplus V^*)^*\equiv V\oplus V^*$. The pair $(V\oplus V^*, \omega)$ is, up to equivalences, the unique \emph{linear symplectic space} of dimension $2n$.
When coordinates are required, we fix a basis $\{e_i\}_{i=1,\dots,n}$ in $V$ and we consider its dual  $\{\epsilon^i\}_{i=1,\dots,n}$ in $V^*$. If the results
do not depend on the ground field, we leave it unspecified---that is, it may be either $\R$ or $\C$.

In the above coordinates, the matrix of $\omega $ is
\begin{equation}\label{eqSympMat}
I:=\left(\begin{array}{cc}0 & -\Id \\\Id & 0\end{array}\right),
\end{equation}
and it is known as the \emph{standard symplectic matrix}. Indeed,
\begin{gather*}
 V\oplus V^* \stackrel{\omega}{\longrightarrow}  V^*\oplus V=( V\oplus V^*)^* ,\\
 v+\alpha \longmapsto -\alpha +v ,
\end{gather*}
with respect to the bases $(e_1,\dots, e_n,\epsilon^1,\dots,\epsilon^n)$ and $(\epsilon^1,\dots,\epsilon^n,e_1,\dots, e_n)$ of $V\oplus V^*$ and $V^*\oplus V$,
respectively.

We stress that \eqref{eqSympMat} is \emph{not} a rotation matrix---it is the matrix corresponding to the two-form
\begin{equation}\label{eqDefFormaSimplettica}
\omega=\epsilon^i\wedge e_i.
\end{equation}
Observe that in \eqref{eqDefFormaSimplettica}, as well as in the rest of the paper, we use the Einstein convention for repeated indexes, unless otherwise
specified.

\subsection{Definition of the Lagrangian Grassmannian}\label{secDefLagGrass}
It is well-known that the set
\begin{equation}\label{eqGrassNonLagr}
\Gr(n,V\oplus V^*):=\{L\subset V\oplus V^*\mid L\textrm{ linear subspace, } \dim L=n\}
\end{equation}
possesses the structure of an $n^2$-dimensional smooth manifold, known as  the \emph{Grassmannian $($manifold/variety\/$)$} (see, e.g.,
\cite[Lecture~6]{harris1992algebraic} for an algebro-geometric proof or \hbox{\cite[Lemma 5.1]{milnor1974characteristic}} for a differential-geometric proof).
The key is the injective map\footnote{Observe that $ \arctan h$ is nothing but the graph of $h$. The symbol ``$\arctan$'' has been chosen in order to be
consistent with  Smith's contribution to this very volume, see \cite[\S 2]{Smith2019}.}
\begin{gather}
V^*\otimes V^*=\Hom(V,V^*) \longrightarrow  \Gr(n,V\oplus V^*) ,\nonumber\\
h \longmapsto \arctan h:=\Span{v+h(v)\mid v\in V},\label{eqArcTanMap}
\end{gather}
which allows one to define an $n^2$-dimensional chart in $\Gr(n,V\oplus V^*)$. This chart is also dense---whence the name \emph{big cell} which we shall
use from now on.

If $h=h_{ij}\epsilon^i\otimes\epsilon^j$, then
\begin{equation}\label{eqArcTanCoord}
\arctan h:=\Span{e_i+h_{ij}\epsilon^j\mid i=1,\dots,n} .
\end{equation}
Let us impose that $\arctan h$ be \emph{isotropic} with respect to the two-form $\omega$, that is
\begin{equation}\label{eqArctTanIsotrop}
\omega|_{\arctan h}\equiv 0.
\end{equation}
In view of \eqref{eqArcTanCoord}, condition \eqref{eqArctTanIsotrop} reads
\begin{equation}
\omega(e_{i_1}+h_{i_1j_1}\epsilon^{j_1},e_{i_2}+h_{i_2j_2}\epsilon^{j_2})=0, \quad\forall i_1,i_2=1,\dots, n.
\end{equation}
Since $\omega(e_{i_1}+h_{i_1j_1}\epsilon^{j_1},e_{i_2}+h_{i_2j_2}\epsilon^{j_2})=h_{i_1i_2}-h_{i_2i_1}$, it is obvious that \eqref{eqArctTanIsotrop} is
fulfilled if and only if the matrix $h_{ij}$ is symmetric, that is $h\in S^2V^*$.

\begin{figure}
\centerline{\epsfig{file=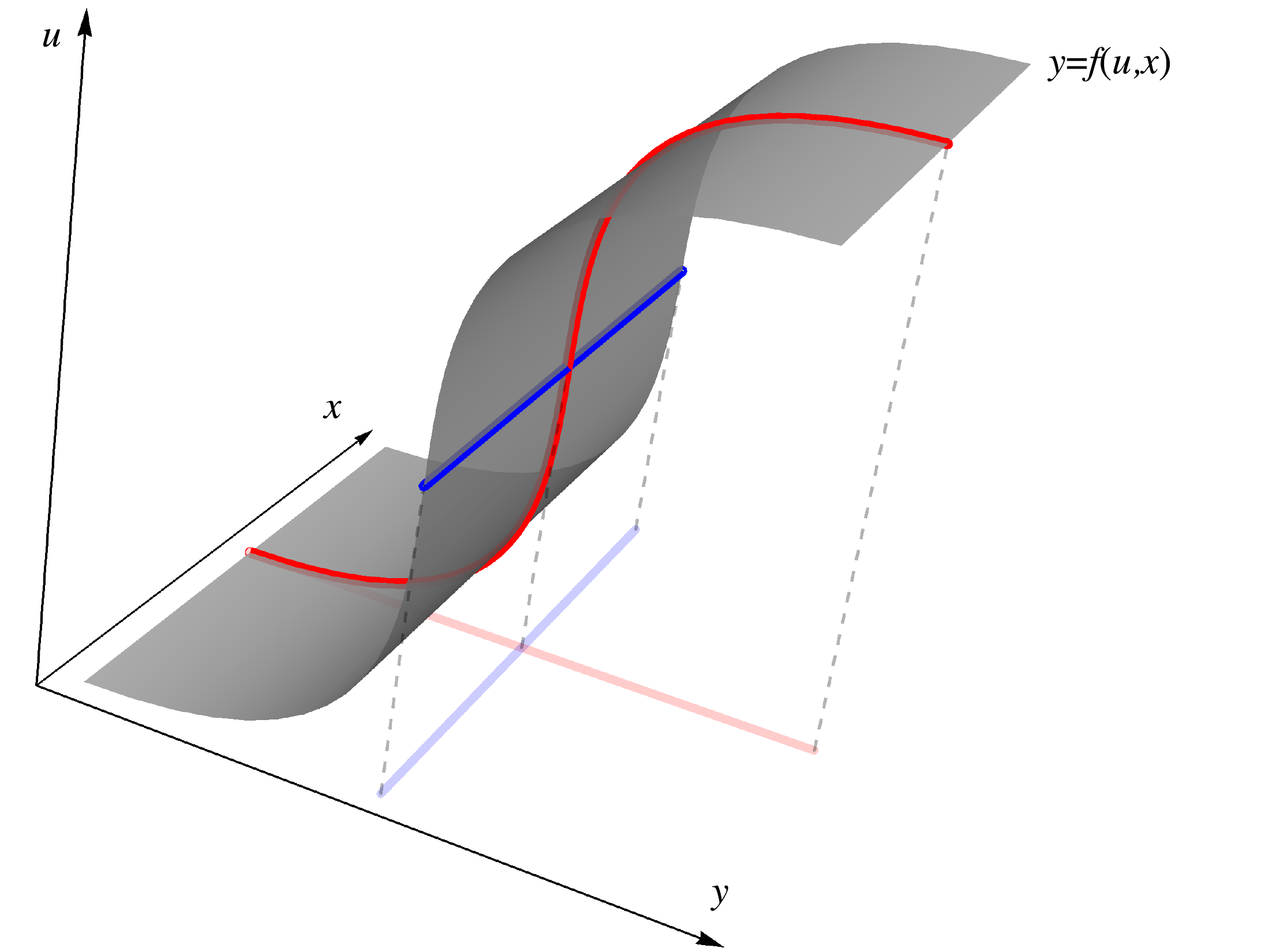,width=0.8\textwidth}}\caption{In the framework of contact geometry, the distinction between dependent ($u$) and independent ($x,y$)  variables simply disappears. The smooth surface depicted here---the graph  of  $y=f(u,x)$, with $f$ smooth---cannot be interpreted as a ``regular'' solution in the sense  of a function $u=u(x,y)$. There is a locus, highlighted as a blue line, where the tangent planes to the surface project degenerately to the $(x,y)$--plane (e.g., the projection of the red curve crosses  the projection of the blue line at zero speed, regardlessly of its parametrisation). Such a locus is called a \emph{singularity} of the solution. We warn the reader that the terminology is misleading, since the surface is perfectly smooth. \label{fig}}
\end{figure}

The \emph{Lagrangian Grassmannian $($manifold/variety\/$)$} $\LL(n,V\oplus V^*)$ can be defined as the closure of the subset $S^2V^*$ of the big cell
$V^*\otimes V^*$. From this point of view, $\LL(n,V\oplus V^*)$ is a \emph{compactification} of the space of symmetric forms on $V$. In the geometric theory
of PDEs based on jet spaces \cite{MR861121}, the additional ``points at infinity'' correspond to the so-called \emph{singularities of solutions}
\cite[\S 2.2]{Vitagliano2013}, see Fig.~\ref{fig}.

\abs0.15 \aabd1
From now on, the open subset
\begin{equation}\label{eqBigCellLG}
 S^2V^*\subset \LL(n,V\oplus V^*)
\end{equation}
will be referred to as the \emph{big cell} of the Lagrangian Grassmannian $\LL(n,V\oplus V^*)$.

\subsection{Coordinate-free definition of the Lagrangian Grassmannian}
On a deeper conceptual level, the symplectic form $\omega$ can be found by decomposing the space of two-forms on $V\oplus V^*$ into  $\GL(V)$-irreducible
representations:
\begin{equation}
\wedge^2(V\oplus V^*)=\wedge^2V\oplus (V \otimes V^*)\oplus \wedge^2V^*.
\end{equation}
Then $\omega$ is precisely the element of the left-hand side corresponding to the element $0+\id+0$ of the right-hand side. Then one can set
\begin{equation}
\LL(V\oplus V^*):=\bigl\{L\subset V\oplus V^*\mid L\textrm{ linear subspace, } \dim L=n,\ \omega|_L\equiv 0\bigr\}.
\end{equation}

\subsection{The Pl\"ucker embedding}
While $S^2V^*$ provides a convenient \emph{local} description of $\LL(V\oplus V^*)$, the rich \emph{global} geometry of $\LL(V\oplus V^*)$ is invisible from
the point of view of the big cell. Global features become evident when the object is embedded into a ``flat'' environment. In the   present case the role
of such an environment is played by an appropriate projective subspace of $\p(\wedge^n(V\oplus V^*))$.

The trick to obtain the desired embedding consists in regarding an $n$-dimensional subspace ${L\in\LL(V\oplus V^*)}$ as a \emph{line} in
$\wedge^n(V\oplus V^*)$. Indeed, a basis $\{l_1,\dots,l_n\} $ of $L$ defines, up to a projective factor, a unique (nonzero) $n$-vector $l_1\wedge\ldots\wedge
l_n$. The projective class of the latter is then unambiguously associated with $L$, and we will call it the \emph{volume} of $L$ and denote it by $ \vol(L)$.

The map
\begin{gather}
\LL(V\oplus V^*) \longrightarrow  \p(\wedge^n(V\oplus V^*)),\label{eqImmPluck}\\
L \longmapsto  \vol(L),\nonumber
\end{gather}
is called the \emph{Pl\"ucker embedding}. The basic properties of \eqref{eqImmPluck} are easily checked.

First, the element $\vol(L)$ is represented by a \emph{decomposable} $n$-vector, that is an $n$-vector $\xi$ satisfying the equation $\xi\wedge\xi=0$.
The latter is a quadratic condition, symmetric for $n$ even and skew-symmetric for $n$ odd.

Second, the representative $\xi$ is transversal to $\omega$, in the sense that
\begin{equation}
\iota_\omega(\xi):=\omega\ins\xi=0.
\end{equation}
This means that, in fact $\xi$  belongs to the linear subspace
\begin{equation}
\ker\iota_\omega=\ker\bigl(  \wedge^n(V\oplus V^*)\stackrel{\iota_\omega}{\longrightarrow} \wedge^{n-2}(V\oplus V^*)  \bigr)
\end{equation}
of $\wedge^n(V\oplus V^*)$.

Third, \eqref{eqImmPluck} is injective.

\subsection{The Pl\"ucker embedding space}

The fact that $\ker\iota_\omega$ is \emph{not} the smallest linear subspace of $\wedge^n(V\oplus V^*)$ whose projectivization contains the image of
\eqref{eqImmPluck}  is less evident and requires more care.

To this end, observe that
\begin{equation}
\wedge^n(V\oplus V^*)=\bigoplus_{i=0}^n\wedge^i(V^*)\otimes \wedge^{n-i}(V)\simeq\bigoplus_{i=0}^n\wedge^i(V^*)\otimes \wedge^i(V^*),
\end{equation}
in view of the Poincar\'e duality $\wedge^{n-i}(V)=\Hom(\wedge^{i}(V),\wedge^{n}(V))$. It is then easy to realise that the representative $\xi$ of $\vol(L)$
belongs to
\begin{equation}\label{eqSpazQuasMin}
\bigoplus_{i=0}^nS^2(\wedge^i(V^*))
\end{equation}
via the map \eqref{eqImmPluck}.
Puzzlingly enough, \eqref{eqSpazQuasMin} is not yet the minimal subspace we were looking for. Though  it is so for $n=2,3$. Let $n=2$ and
let $L$ be the Lagrangian 2-plane \eqref{eqArcTanCoord} corresponding to the symmetric $2\times 2$ matrix $(h_{ij})$. Then $\vol(L)=[\xi]$, with
\begin{equation}
\xi=e_1\wedge e_2+h_{11}e_1\wedge\epsilon^1+h_{12}e_1\wedge\epsilon^2+h_{21}e_2\wedge\epsilon^1+h_{22}e_2\wedge\epsilon^2+\det(h_{ij})\epsilon^1\wedge\epsilon^2
\end{equation}
as an element of $\wedge^2(V\oplus V^*)$, and
\begin{equation}
\xi=1+h_{ij}\epsilon^i\odot\epsilon^j+\det(h_{ij})\epsilon^1\wedge\epsilon^2
\end{equation}
as an element of \eqref{eqSpazQuasMin}. That is,
\begin{equation}\label{eqMinorsOrd2}
(1,h_{11}, h_{12},h_{22},h_{11}h_{22}-h_{12}^2)
\end{equation}
are the coordinates of $\xi$ in the standard basis of  \eqref{eqSpazQuasMin}. Observe that in \eqref{eqMinorsOrd2} there appear all the minors  of the matrix
$(h_{ij})$, namely: minors of order 0 (the constant 1), minors of order 1 (the very entries of the matrix) and minors of order 2 (the determinant).

Similarly, for $n=3$, one finds
\begin{equation}
\xi=1+h_{ij}\epsilon^i\odot\epsilon^j+h_{ij}^\#(\epsilon^1\wedge\widehat{\epsilon^i}\wedge\epsilon^3)\odot (\epsilon^1\wedge\widehat{\epsilon^j}\wedge\epsilon^3) +\det(h_{ij})\epsilon^1\wedge\epsilon^2\wedge\epsilon^3,
\end{equation}
where $(h_{ij}^\#)$ denotes the cofactor matrix of $(h_{ij})$ and the hat indicates a removed element. One gets then the 14 coordinates
\begin{equation}\label{eqMinorsOrd3}
(1, h_{11}, h_{12}, h_{13}, h_{22}, h_{23}, h_{33}, h_{11}^\#, h_{12}^\#, h_{13}^\#, h_{22}^\#, h_{23}^\#, h_{33}^\#,\det(h_{ij}))
\end{equation}
corresponding to the point $\xi$.

It is then easy to realise that, in the case $n=4$, one has
\begin{equation}\label{eqPluckCoord4x4}
(1,\dots,h_{ij},\dots,2\times2\ \textrm{minors} ,\dots, h_{ij}^\#,\dots, \det(p_{ij})).
\end{equation}
A $4\times 4$ symmetric matrix contains exactly: 1 minor of order 0 and of order 4, $\frac{4\cdot 5}{2}=10$ minors of order 1 and 3, $\frac{6\cdot 7}{2}=21$
minors of order 2, where $6={4\choose 2}$ is the number of choices of 2 rows (columns). Therefore, \eqref{eqPluckCoord4x4} consists  exactly of
$1+10+21+10+1=43$ entries. The subtle point here is that, as opposed to the cases $n=2$ and $n=3$, \emph{not all the minors are linearly independent}.
More precisely there is exactly one linear combination of $2\times 2$ minors, namely
\begin{equation}\label{eqRelazioneMinori4x4}
-(h_{13}h_{24}-h_{14}h_{23})+(h_{12}h_{34}-h_{14}h_{23})-(h_{12}h_{34}-h_{13}h_{24}),
\end{equation}
which vanishes \cite{MO209058}. Therefore, there is a (proper, for $n\geq 4$) linear subspace of \eqref{eqSpazQuasMin}, henceforth denoted by
\begin{equation}\label{eqSpazioMinimalePlucker}
\wedge^n_0(V\oplus V^*):=\bigoplus_{i=0}^nS_0^2(\wedge^i(V^*)),
\end{equation}
which contains all the $\xi$'s and it is minimal with respect to this property. Summing up,
\begin{align*}
n=2 &\Rightarrow  \wedge^2_0(V\oplus V^*)=\wedge^0(V^*)\oplus S^2(V^*)\oplus \wedge^2(V^*)\textrm{ has dimension }5,\\
n=3 &\Rightarrow  \wedge^3_0(V\oplus V^*)=\wedge^0(V^*)\oplus S^2(V^*)\oplus S^2(\wedge^2(V^*))\oplus \wedge^3(V^*) \textrm{ has dimension }14,\\
n=4 &\Rightarrow  \wedge^4_0(V\oplus V^*)=\wedge^0(V^*)\oplus S^2(V^*)\oplus S^2_0(\wedge^2(V^*))\oplus S^2(\wedge^3(V^*))\oplus \wedge^4(V^*)\\
&\hphantom{\Rightarrow}\textrm{ has dimension }42.
\end{align*}
Therefore, the minimal projective embedding of the 3-(resp., 6- and 10-)dimensional Lagrangian Grassmannian $\LL(2,4)$ (resp., $\LL(3,6)$ and $\LL(4,8)$)
is $\p^4$ (resp., $\p^{13}$ and $\p^{41}$).

In general, to find the Pl\"ucker embedding space of the $\frac{n(n+1)}{2}$-dimensional Lagrangian Grassmannian $\LL(n,2n)$, one has to count how many minors
a symmetric $n\times n$ matrix possesses, minors of order 0 and $n$ included. This, in principle, is an easy task. The problem is to look for dependencies
of the form \eqref{eqRelazioneMinori4x4} among minors. The number of minors needs to be diminished by the number of these relations. The result, further
decreased by one, represents the (projective) dimension of the sought-for space. In Section \ref{SecRepStuff} below we explain how the very same space can be
obtained by exploiting the theory of representations.

\subsection{The Pl\"ucker relations}
Expressions \eqref{eqMinorsOrd2}, \eqref{eqMinorsOrd3} and \eqref{eqPluckCoord4x4} represent the parametric description of $\LL(2,4)$, $\LL(3,6)$ and
$\LL(4,8)$ in $\p^4$, $\p^{13}$ and $\p^{41}$, respectively. Let us denote by
\begin{equation}
[z_0:z_1:\ldots :z_N]
\end{equation}
the standard projective coordinates on $\p^N$. Then it is easy to realise that points of $\LL(2,4)$ satisfy the quadratic relation
\begin{equation}\label{eqL24}
z_1z_3-z_2^2-z_0z_4=0,
\end{equation}
capturing the fact that, on $\LL(2,4)$, the fourth coordinate is the determinant of the symmetric matrix whose entries are $z_1$, $z_2$ and $z_3$. Observe that
$z_0$, which is 1 on the big cell of $\LL(2,4)$, serves the sole purpose of homogenising the relation
\begin{equation}
z_4=\det\left(\begin{array}{cc}z_1 & z_2 \\z_2 & z_3\end{array}\right).
\end{equation}
Indeed, \eqref{eqL24} is the equation cutting out $\LL(2,4)$ in $\p^4$. Therefore, $\LL(2,4)$ is a quadric.

\abs-0.1 \aabd-1
The codimension of $\LL(3,6)$ in $\p^{13}$, on the contrary, is quite high: 7. The corresponding equations are essentially the Laplace rule for the determinant
of a symmetric $3\times 3$ matrix. More precisely, if
\begin{equation}
Z:=\left(\begin{array}{ccc}z_1 & z_2 & z_3 \\z_2 & z_4 & z_5 \\z_3 & z_5 & z_6\end{array}\right),
\end{equation}
then the seven quadratic equations
\begin{align*}
z_0z_{7} &= Z^\#_{11}\, \\
z_0z_{8} &= Z^\#_{12}\, \\
z_0z_{9} &= Z^\#_{13}\, \\
z_0z_{10} &= Z^\#_{22}\, \\
z_0z_{11} &= Z^\#_{23}\, \\
z_0z_{12} &= Z^\#_{33}\, \\
 z_0z_{13} &= \tr(Z\cdot Z^\#)
\end{align*}
cut out $\LL(3,6)$ in $\p^{13}$. Similarly, it can be proved that $\LL(n,2n)$ is cut out by quadratic relations in its own Pl\"ucker embedding space
$\p(\wedge_0^n(V\oplus V^*))$, similar to \eqref{eqL24} and the seven equations above (see, e.g., \cite[Theorem 14.6]{miller2004combinatorial}). These
relations are usually referred to as the \emph{Pl\"ucker relations}, whereas the expressions  \eqref{eqMinorsOrd2}, \eqref{eqMinorsOrd3},
\eqref{eqPluckCoord4x4}, as well as the analogous ones for higher values of $n$, are called the \emph{Pl\"ucker coordinates} of the point $L\in\LL(n,2n)$.

\subsection{The dual variety}\label{secDualVAr}
The case of $\LL(2,4)$ is somewhat special in that the dimension of the Pl\"ucker embedding space $\p^4$ exceeds only by one  the dimension of $\LL(2,4)$.
Then the tangent spaces to $\LL(2,4)$ are projective hyperplanes in $\p^4$. The latter form a set, usually denoted by $\p^{4\,\ast}$ and called the \emph{dual}
of $\p^4$, which is (non-canonically) identified with $\p^4$ itself. Thus, the set of tangent hyperplanes to $\LL(2,4)$ constitutes a subset
\begin{equation}
\LL(2,4)^*\subseteq \p^{4\,\ast}
\end{equation}
of the set of all hyperplanes, accordingly called the \emph{dual variety} of $\LL(2,4)$.

We warn the reader that an element $\pi\in \LL(2,4)^*$ is a linearly embedded $\p^3$, which has contact of order one \emph{in some point} with $\LL(2,4)$.
In the language of jets,
\begin{equation}
\LL(2,4)^*=\{ \pi\in\p^{4\,\ast}\mid j^1_x(\pi)=j^1_x(\LL(2,4))\textrm{ in some point }x\in \LL(2,4)\}.
\end{equation}
Therefore, the same $\pi$, which is tangent at $x\in\LL(2,4)$, may intersect unpredictably $\LL(2,4)$ someplace else. In fact, as we shall see later on,
elements of $ \LL(2,4)^*$ allow constructing special hypersurfaces in $\LL(2,4)$ called \emph{hyperplane sections}.  Hence, the notion of an element of the
dual variety is different from (though related to) the notion of a tangent space to $\LL(2,4)$, that is a fibre of an abstract linear bundle of rank 3.

Another peculiarity of $\LL(2,4)$ is that---exceptionally among all Lagrangian Grass\-mann\-ians---its dual variety is smooth and canonically isomorphic
to $\LL(2,4)$ itself. That is, it is cut by the very same equation \eqref{eqL24}, appropriately interpreted as an equation in $\p^{4\,\ast}$. This also follows
from the fact that $\LL(2,4)$ has codimension one in $\p^4$, and that $\LL(2,4)$ is smooth. Indeed, to any point $x\in\LL(2,4)$ one associates the
\emph{unique} element $\pi_x\in \p^{4\,\ast}$ which is tangent to $\LL(2,4)$ at $x$. This realises the desired one-to-one correspondence.

The case of $\LL(3,6)$ is already much more involved. Indeed, at any $x\in\LL(3,6)$ there is certainly a unique tangent \emph{$6$-dimensional subspace} but there
does not need to be a unique tangent \emph{hyperplane} (i.e., a 12-dimensional subspace). Actually, there is a 6-dimensional family of them, making the dual
variety $\LL(3,6)^*$ 12-dimensional. It can be proved that it is cut out by a single quartic relation in $\p^{13\,\ast}$ \cite[Section 5]{Russo2019}.

\abs0.1 \aabd1
There is still a certain correspondence between $\LL(3,6)$ and its dual $\LL(3,6)^*$. The former is isomorphic to the singular locus\footnote{By
\emph{singular locus} of an algebraic variety $X=\{f_1=0,\dots, f_m=0\}$ of codimension $m$ we mean the subset of $X$ where the differentials $df_1,\dots,
df_m$ are not linearly independent.} of the latter, viz.
\begin{equation}\label{eqL36isoDualSinglL36}
\LL(3,6)\equiv\Sing(\LL(3,6)^*).
\end{equation}
This should help  to convince oneself of the validity of \eqref{eqL36isoDualSinglL36}. Let us describe an element $\pi\in\p^{14\,\ast}$ by projective
coordinates
\begin{equation}
\pi\equiv [A:\ldots :B^{ij}:\ldots :\ldots: C^{ij}:\ldots: D].
\end{equation}
Then, the intersection $\pi\cap\LL(3,6)$, in the Pl\"ucker coordinates  \eqref{eqMinorsOrd3}, is given by
\begin{equation}\label{eqProtoHypSect}
A+B^{ij}h_{ij}+C^{ij}h^\#_{ij}+D\det(h_{ij})=0.
\end{equation}
The key remark is that a particular case of an expression of the form \eqref{eqProtoHypSect} can be obtained by means of \emph{another} symmetric
$3\times 3$ matrix, say $H=(H_{ij})$. More precisely,
\begin{equation}\label{eqProtoGoursatMAEQ}
\det(h-H)=0
\end{equation}
is a particular form of equation \eqref{eqProtoHypSect} above, where the \emph{fourteen} coefficients $A, B^{ij}, C^{ij},D$ depend on the \emph{six}
coefficients $H_{ij}$. It is not hard to realise that, after the substitutions
\begin{equation}\label{eqCoefficientiDuali}
A:=\det H,\quad B:=H^\#,\quad C:=H,\quad D:=1,
\end{equation}
the equation \eqref{eqProtoHypSect} becomes \eqref{eqProtoGoursatMAEQ}. On the top of that, the hyperplane $\pi_H$, with coefficients given by
\eqref{eqCoefficientiDuali} is tangent to $\LL(3,6)$. This is not hard to see: the left-hand side of \eqref{eqProtoGoursatMAEQ}, regarded as a function of $h$,
vanishes at $h=H$, \emph{together with its first derivatives}.

In other words $\pi\in \LL(3,6)^*$ because $\pi_H$ is tangent to $\LL(3,6)$ at the point $x_H$ given, in the  coordinates  \eqref{eqMinorsOrd3},  by $H$ itself.
The correspondence
\begin{equation}\label{eqCooDuality}
\pi_H\longmapsto x_H
\end{equation}
basically allows us to regard the same matrix $H$ as a special element $\pi_H$ of $\LL(3,6)^*$ as well as an element $x_H$ of $\LL(3,6)$ itself, thus realising
the desired isomorphism \eqref{eqL36isoDualSinglL36}.

The duality \eqref{eqCooDuality} manifests itself for any $\LL(n,2n)$, though the isomorphism \eqref{eqL36isoDualSinglL36} now must be recast as
\begin{equation}
\LL(n,2n)\equiv\underbrace{ \Sing(\cdots\Sing(}_{n-2\textrm{ times}}\LL(n,2n)^* )\cdots ).
\end{equation}
The underlying structure responsible of this duality is the natural bilinear form
\begin{gather}
\wedge_0^n(V\oplus V^*) \times \wedge_0^n(V\oplus V^*) \longrightarrow \wedge^{2n}(V\oplus V^*),\nonumber\\
(\alpha,\beta) \longmapsto \alpha\wedge\beta,\label{eqFormSimmPluckSpace}
\end{gather}
which is scalar-valued and is symmetric (resp., skew-symmetric) for $n$ even (resp., odd). Indeed, the above-defined bilinear form is non-degenerate, thus
allowing a point-hyper\-plane correspondence in the (de-projectivised) Pl\"ucker embedding space $\wedge_0^n(V\oplus V^*) $.  After projectivisation, this
correspondence coincides precisely with \eqref{eqCooDuality}. It is interesting to notice that such a correspondence is equivalent to the fact that the cone
over $\LL(n,2n)$ be \emph{isotropic} with respect to \eqref{eqFormSimmPluckSpace}.

The dual $\p(\wedge_0^n(V\oplus V^*) )^*$ of the Pl\"ucker embedding space parametrises the hyperplane sections of the Lagrangian Grassmannian, which
correspond to the so-called \emph{Monge--Amp\`ere equations}, see Section \ref{secHypMAEs} below. The  stratification of $\p(\wedge_0^n(V\oplus V^*) )^*$
by the dual variety $\LL(n,2n)^*$ and its singular loci will correspond to special ($\Sp_{2n}$-invariant, see next Section \ref{secNatGrAct}) classes
of such PDEs (see Section \ref{secLowDimExamples}).

The study of these special classes of PDEs corresponds precisely to the study of the orbits in $\p(\wedge_0^n(V\oplus V^*) )^*$  of the natural groups acting
on $\wedge_0^n(V\oplus V^*)$ (see Section \ref{secNatGrAct} below). In particular, there is a unique close orbit with respect to the symplectic group, and
this is precisely the ``very singular'' locus $\LL(n,2n)$ (see Section \ref{secHomStructLG}).

\subsection{Natural group actions on $\LL(V\oplus V^*)$}\label{secNatGrAct}
Recall that the ``arctangent map'' \eqref{eqArcTanMap} allowed us to define a canonical embedding of $S^2V^*$ into the Lagrangian Grassmannian
$\LL(V\oplus V^*)$, whose image corresponds to  the {big cell} \eqref{eqBigCellLG}. Let us further restrict our scope by considering only the
\emph{non-degenerate elements} of the big cell. That is, the open subset
\begin{equation*}
\cU:=\{ h\in S^2V^*\mid \det(h)\neq 0\}\subset \LL(V\oplus V^*).
\end{equation*}
One obvious group action  on $\cU$ is easily found. Indeed, any element $D\in \GL(V)$ acts naturally on  symmetric forms,
\begin{equation}\label{eqNatActGLV}
h \longmapsto D^t\cdot h\cdot D,
\end{equation}
where the same symbol $h$ denotes both the matrix and the form itself. As a matter of fact, \eqref{eqNatActGLV} acts on the whole big cell, preserving
the subset $\cU$.

Another group action on $\cU$ is due to the linear structure of the big cell. Indeed, an element $C$ of the \emph{Abelian group} $ S^2V^*$ can act on $S^2V^*$
itself as a \emph{translation}:
\begin{equation}\label{eqActTransS2VStar}
h \longmapsto h+ C.
\end{equation}
Observe that, unlike \eqref{eqNatActGLV}, the action \eqref{eqActTransS2VStar} does not preserves $\cU$.

One last, somewhat less evident, group action on $\cU$ is given by
\begin{equation}\label{eqActTransS2V}
h \longmapsto h\cdot(\Id+B\cdot h)^{-1},
\end{equation}
where now $B\in S^2V$. Above correspondence \eqref{eqActTransS2V} can be explained as follows. There is an analogue of the arctangent map \eqref{eqArcTanMap},
defined on $S^2V$, instead of $S^2V^*$, that is
\begin{equation}
S^2V\ni k\longmapsto \underline{\arctan}(k):=\Span{\alpha+k(\alpha)\mid \alpha\in V^*}.
\end{equation}
Observe that the common image of $\arctan$ and $\underline{\arctan}$ is precisely $\cU$. Therefore, \eqref{eqActTransS2V} is nothing but the translation
\begin{equation}\label{eqActTransS2Vdemisitfied}
h^{-1} \longmapsto h^{-1} +B ,
\end{equation}
by $B$ of $h^{-1}$, understood as an element of $S^2V$ via $\underline{\arctan}$. Indeed,
\begin{align*}
\lefteqn{\kern-30pt \underline{\arctan}(h^{-1} +B ) = \Span{\alpha+(h^{-1} +B)(\alpha)\mid \alpha\in V^*}}\\
&=  \Span{h(v)+(h^{-1} +B)(h(v))\mid v\in V}
=  \Span{h(v)+(\Id +B\cdot h)(v)\mid v\in V}\\
&=  \Span{h((\Id +B\cdot h)^{-1}(v))+(\Id +B\cdot h)((\Id +B\cdot h)^{-1}(v))\mid v\in V}\\
&=  \Span{v+(h(\Id +B\cdot h)^{-1})(v)\mid v\in V}
=\arctan(h\cdot (\Id +B\cdot h)^{-1}),
\end{align*}
where we used the facts that $h$ is invertible and that, at least locally around 0, $\Id +B\cdot h$ is invertible as well.

\abs-0.1 \aabd-1
The three actions \eqref{eqNatActGLV}, \eqref{eqActTransS2VStar} and \eqref{eqActTransS2V} above may seem accidental and unrelated. On the contrary, they
share a common background. Consider a linear transformation of $V\oplus V^*$, represented, in the aforementioned basis, by the $(2n)\times(2n)$ matrix
\begin{equation}
M:=\left(\begin{array}{cc}A & B \\C & D\end{array}\right).
\end{equation}
In the same basis, the $(2n)\times n$ matrix
\begin{equation}\label{eqMatrixArcTanH}
\left(\begin{array}{c}\Id_n \\h\end{array}\right)\equiv\arctan(h)
\end{equation}
represents the $n$-dimensional linear subspace $\arctan(h)$. Indeed, the $n$ columns of the matrix \eqref{eqMatrixArcTanH} corresponds to the $n$ generators
appearing in the definition \eqref{eqArcTanCoord} of $\arctan(h)$. Observe that
\begin{equation}
\left(\begin{array}{c}\Id_n \\h\end{array}\right), \quad \left(\begin{array}{c}A \\h\cdot A\end{array}\right)
\end{equation}
represent the same subspace, for any $A\in\GL(V)$.

We need now to make  the crucial assumption that $M$ belongs to a small neighbourhood of the identity. This allows us to act by  $M$ on $\arctan(h)$ as follows:
\begin{align}
M\cdot \arctan(h)&=\left(\begin{array}{cc}A & B \\C & D\end{array}\right)\cdot \left(\begin{array}{c}\Id_n \\h\end{array}\right) \nonumber\\
&=\left(\begin{array}{c}A+B\cdot h \\ C+D\cdot h\end{array}\right)
=\left(\begin{array}{c}\Id_n \\ (C+D\cdot h)(A+B\cdot h)^{-1}\end{array}\right) \nonumber\\
&=  \arctan((C+D\cdot h)(A+B\cdot h)^{-1}).\label{eqAzioneMSuArctan}
\end{align}
Directly from the definition of a Lagrangian subspace of $V\oplus V^*$ it follows that   $(C+D\cdot h)(A+B\cdot h)^{-1}$ is again a symmetric form
if and only if the transformation $M$ preserves the symplectic form $\omega$, that is,
\begin{equation}\label{eqSimplTransfMatrix}
M^t\cdot I\cdot M=I,
\end{equation}
where $I$ is the \emph{symplectic matrix} \eqref{eqSympMat}.
A matrix $M$ fulfilling \eqref{eqSimplTransfMatrix} is called \emph{symplectic transformation}. The three matrices
\begin{equation}\label{eqTreMAtriciSimplettiche}
\left(\begin{array}{cc}(D^t)^{-1} & 0 \\0 & D\end{array}\right),\quad \left(\begin{array}{cc}1 & 0 \\C & 1\end{array}\right),\quad
\left(\begin{array}{cc}1 & B \\0 & 1\end{array}\right),
\end{equation}
with $D\in\GL(V)$, $C\in S^2V^*$ and $B\in S^2V$
are easily checked to satisfy \eqref{eqSimplTransfMatrix} and they correspond to the actions \eqref{eqNatActGLV}, \eqref{eqActTransS2VStar} and
\eqref{eqActTransS2V}, respectively.

Actually, the three matrices \eqref{eqTreMAtriciSimplettiche} generate the entire subgroup
\begin{equation}
\Sp_{2n} \equiv \Sp(V\oplus V^*)\subset\GL(V\oplus V^*)
\end{equation}
of symplectic transformations, that is what is usually called the \emph{symplectic group}. Such an ``inner structure'' of the symplectic group becomes
even more evident on the level of the corresponding   Lie algebras, viz.
\begin{equation}\label{eqGradSpEnneAlgebra}
\sp_{2n}\equiv\sp(V\oplus V^*)=S^2V^*\oplus\gl(V)\oplus S^2V.
\end{equation}
This structure is the source of all the structures on $\LL(V\oplus V^*)$ we shall find later on. The homogeneous one, to begin with.

\abs0.1 \aabd1
\subsection{The homogeneous structure of $\LL(V\oplus V^*)$}\label{secHomStructLG}

Formula \eqref{eqAzioneMSuArctan} immediately shows that
\begin{equation}
M\cdot\arctan(0)=\arctan(0)\Leftrightarrow A=\Id_n,\ C=0.
\end{equation}
In other words, the stabiliser subgroup
\begin{equation}\label{eqPrimaDefPi}
P:=\Stab_{\Sp_{2n}}(\arctan(0))=\left\{  \left(\begin{array}{cc}\Id_n & B \\0 & D\end{array}\right)\Bigm| B\in S^2V,\ D\in\GL(V)  \right\},
\end{equation}
which coincides with the semidirect product
\begin{equation}\label{eqDefSubGroupPi}
P=\GL(V)\rtimes S^2V,
\end{equation}
encompasses the transformations of the form \eqref{eqNatActGLV} and \eqref{eqActTransS2V}. Those of the form \eqref{eqActTransS2VStar}, that is $S^2V^*$
acting by translations on the big cell $S^2V^*$, clearly allow us to move the origin  $\arctan(0)$ into any other point $\arctan(h)$ of the big cell. So,
the orbit $\Sp_{2n}/P$ contains the big cell $S^2V^*$. However, since $\Sp_{2n}$ is compact and $P$ is closed,   it must be
\begin{equation}
\Sp_{2n}/P=\overline{S^2V^*}=\LL(n,2n).
\end{equation}
That is, the $\Sp_{2n}$-action is transitive and $\LL(n,2n)$ is a homogeneous space of the Lie group $\Sp_{2n}$.

Having ascertained the transitivity of the $\Sp_{2n}$-action, we can switch to the local point of view and analyse the infinitesimal action of $\sp_{2n}$.
Assume that
\begin{equation}
M_\epsilon=\left(\begin{array}{cc}A_\epsilon & B_\epsilon \\C_\epsilon & D_\epsilon\end{array}\right)
\end{equation}
passes through the identity at $\epsilon=0$, and differentiate   formula \eqref{eqAzioneMSuArctan}:
\begin{align}
\lefteqn{\frac{\dd}{\dd \epsilon}\biggr|_{\epsilon=0}M_\epsilon\cdot h =
\frac{\dd}{\dd \epsilon}\biggr|_{\epsilon=0}(C_\epsilon+D_\epsilon\cdot h)(A_\epsilon+B_\epsilon\cdot h)^{-1}}\nonumber\\
&=(\dot{C}_0+\dot{D}_0\cdot h)(A_0+B_0\cdot h)^{-1}-(C_0+D_0\cdot h)(\dot{A}_0+\dot{B}_0\cdot h)(A_0+B_0\cdot h)^{-2}\nonumber\\
&=\dot{C}_0+\dot{D}_0\cdot h-h\cdot (\dot{A}_0+\dot{B}_0\cdot h),\label{eqInfinitesimalActionOfEmme}
\end{align}
as well as \eqref{eqSimplTransfMatrix}:
\begin{align}
0&=\frac{\dd}{\dd \epsilon}\biggr|_{\epsilon=0} M^t_\epsilon\cdot I\cdot M_\epsilon=\dot{M}^t_0\cdot I\cdot M_0+M^t_0\cdot I\cdot \dot{M}_0\nonumber\\
&=\dot{M}^t_0\cdot I+I\cdot \dot{M}_0\nonumber\\
&=\left(\begin{array}{cc}\dot{A}^t_0 & \dot{C}_0^t\\ \dot{B}_0^t & \dot{D}_0^t\end{array}\right)\cdot \left(\begin{array}{cc}0 & -\Id \\
\Id & 0\end{array}\right)+\left(\begin{array}{cc}0 & -\Id \\\Id & 0\end{array}\right)\cdot \left(\begin{array}{cc}\dot{A}_0 & \dot{B}_0\\
\dot{C}_0 & \dot{D}_0\end{array}\right)\nonumber\\
&= \left(\begin{array}{cc}\dot{C}^t_0 & -\dot{A}_0^t\\ \dot{D}_0^t & -\dot{B}_0^t\end{array}\right)+\left(\begin{array}{cc}-\dot{C}_0 & -\dot{D}_0\\
\dot{A}_0 & \dot{B}_0\end{array}\right).\label{eqInfinitesimalSymplecticMatrix}
\end{align}
From \eqref{eqInfinitesimalSymplecticMatrix} we obtain
\begin{equation}
\sp_{2n}=\left\{\dot{M}_0=\left(\begin{array}{cc}-\dot{D}^t_0 & \dot{B}_0\\ \dot{C}_0 & \dot{D}_0\end{array}\right)\Bigm| \dot{B}_0\in S^2V,\
\dot{C}_0\in S^2V^*,\ \dot{D}_0\in\gl(V)\right\},\label{eqDescrSpEnneAlgebra}
\end{equation}
whence \eqref{eqInfinitesimalActionOfEmme} become
\begin{equation}\label{eqAzioneInfinitesimaleMPuntoZero}
\dot{M}_0\cdot h = \dot{C}_0 +\dot{D}_0\cdot h+h\cdot \dot{D}_0^t-h\cdot \dot{B}_0\cdot h.
\end{equation}
The decomposition \eqref{eqGradSpEnneAlgebra} is implicitly written already in \eqref{eqDescrSpEnneAlgebra}, and it should be interpreted as
a $|1|$-grading,\footnote{See \cite[Definition 3.1.2]{MR2532439} for the general definition of a $|k|$-grading.} i.e.,
\begin{equation}\label{eqGradSpEnneAlgebraBIS}
\sp_{2n}=\underbrace{S^2V^*}_{\deg=-1}\oplus\underbrace{\gl(V)}_{\deg=0}\oplus \underbrace{S^2V}_{\deg=+1}.
\end{equation}
In particular, it follows from \eqref{eqGradSpEnneAlgebraBIS} that both $S^2V^*$ and $S^2V$ are Abelian Lie algebras with a (natural) structure of
$\gl(V)$-module. Accordingly, the subgroup $P$ defined by \eqref{eqDefSubGroupPi} corresponds infinitesimally to the non-negative part of the grading:
\begin{equation}
\gp= \gl(V)\oplus S^2V.
\end{equation}
The remaining part, that is $S^2V^*$, is canonically identified with the tangent space at $\arctan(0)$ to $\LL(n,2n)$:
\begin{equation}\label{eqPrimaIdentificazioneSpazioTangenteInUnPunto}
T_{\arctan(0)}\LL(n,2n)\equiv S^2V^*.
\end{equation}
Identification \eqref{eqPrimaIdentificazioneSpazioTangenteInUnPunto} will play a crucial role in the sequel. We stress here that, due to the presence of
a quadratic term in $h$ in \eqref{eqAzioneInfinitesimaleMPuntoZero}, the isotropy action of $P$ on the tangent space $T_{\arctan(0)}\LL(n,2n)$ reduces
to the natural action of its 0-graded part, that is $\gl(V)$, on $S^2V^*$. The action of its  1-graded part, that is $S^2V$, becomes visible only when
the \emph{principal bundle}
\begin{equation}\label{eqPrincBundleSpn}
\gathered
\xymatrix{
\Sp_{2n}\ar[d]^P\\\LL(n,2n)
}\endgathered
\end{equation}
is identified with a sub-bundle of the \emph{second-order frame bundle} of $\LL(n,2n)$ (see Section \ref{sec2NdFrBund} below).

\subsection{The tautological and the tangent bundle of $\LL(n,2n)$}
We discuss now two important linear bundles that  can be naturally associated with $\LL(n,2n)$---the \emph{tautological} (rank-$n$) bundle and the tangent
bundle (whose rank is $\frac{n(n+1)}{2}$). The key observation is that the latter can be identified with the second symmetric power of the dual of the former.
Definitions can be easily given in terms of associated bundles to the $P$-principal bundle \eqref{eqPrincBundleSpn} introduced above. The key identification,
on the other hand, is more evident from a local perspective.

From the $P$-principal bundle \eqref{eqPrincBundleSpn} one immediately obtains the \emph{linear} bundle   $\Sp_{2n}\times^PS^2V^*$
by letting $P$ act on $S^2V^*$ naturally through its 0-graded part and trivially through the rest. By definition, the associated bundle is precisely
the tangent bundle to $\LL(n,2n)$, viz.
\begin{equation}\label{eqIdentTanS2Associato}
T\LL(n,2n)= \Sp_{2n}\times^PS^2V^*.
\end{equation}
We can regard \eqref{eqIdentTanS2Associato} as a generalisation of \eqref{eqPrimaIdentificazioneSpazioTangenteInUnPunto} above, in the sense that the former,
evaluated at $\arctan(0)$, gives the latter. The very\vadjust{\eject} identification \eqref{eqIdentTanS2Associato}  indicates  also how to define a linear \emph{rank-$n$}
bundle, whose  dualised symmetric square coincides with  the tangent bundle to $\LL(n,2n)$. It suffices to rewrite  \eqref{eqIdentTanS2Associato} as
\begin{equation}\label{eqIdentFundPRE}
T\LL(n,2n)= S^2(\Sp_{2n}\times^P V)^*.
\end{equation}
Indeed, at the right-hand side of \eqref{eqIdentFundPRE}, we see now the symmetric square of the dual of the following rank-$n$ bundle
\begin{equation}\label{eqTautBundPRE}
\gathered
\xymatrix{
L:=\Sp_{2n}\times^P V\ar[d]\\
\LL(n,2n).
}\endgathered
\end{equation}
We call  \eqref{eqTautBundPRE} the \emph{tautological bundle} and we denote it by the symbol $L$. The choice of the letter $L$ is not accidental:
if the same symbol $L$ denotes both the total space of the bundle \eqref{eqTautBundPRE}  over $\LL(n,2n)$ \emph{and} a point $L\in \LL(n,2n)$, then
\begin{equation}\label{eqTautBund}
L_L=L,
\end{equation}
that is, the fibre of $L$ at $L$ is again $L$---whence the modifier ``{tautological}''.
With this notation, \eqref{eqIdentFundPRE} becomes simply
\begin{equation}\label{eqIdentFund}
T\LL(n,2n)=S^2L^*.
\end{equation}
Observe that by
\begin{equation}\label{eqIdentFundPointEll}
T_L\LL(n,2n)=S^2L^*
\end{equation}
we mean  that the bundle identification \eqref{eqIdentFund} has been evaluated at the particular point $L\in\LL(n,2n)$, thus becoming an identification
of linear spaces. The reader should be aware that \eqref{eqIdentFund} is an identification of bundles, whereas  \eqref{eqIdentFundPointEll} is an
identification of linear spaces, in spite of the usage of the same symbol $L$.

The importance of \eqref{eqIdentFund} is that it allows us to speak about the \emph{rank} of a tangent vector to $\LL(n,2n)$, which is the rank of the
corresponding bilinear form on the tautological bundle. In particular, rank-one vectors will be tightly connected to the key notion of a~\emph{characteristic}
of a second-order PDE (see Sections \ref{secRO} and \ref{SecChar}).

For the reader feeling uncomfortable with the language of induced bundles we propose another explanation of the identification \eqref{eqIdentFundPointEll}.
Regard $L$ as a point of the Grassmannian $\Gr(n,V\oplus V^*)$ of $n$-dimensional subspaces of $V\oplus V^*$ (see \eqref{eqGrassNonLagr}) and observe that
the $\arctan$ map \eqref{eqArcTanMap} can be generalised by choosing an arbitrary complement $L^c$ of $L$ in $V\oplus V^*$ and by defining
\begin{equation}\label{eqArcTanElle}
\arctan_L:\Hom(L,L^c)\longrightarrow \Gr(n,V\oplus V^*)
\end{equation}
exactly the same way as $\arctan$. Now the symplectic form $\omega$ allows us to identify $L^c$ with the dual $L^*$, whence $\Hom(L,L^c)$ with $L^*\otimes L^*$.
The differential at 0 of \eqref{eqArcTanElle} gives then an isomorphism between $L^*\otimes L^*$ and $T_L \Gr(n,V\oplus V^*)$, which one shows not to depend
upon the choice of $L^c$. Finally, by similar  reasonings as those in Section \ref{secDefLagGrass}  one  finds out that the subspace $S^2L^*$ corresponds
precisely to the subspace  $T_L\LL(n,2n)$, thus obtaining \eqref{eqIdentFundPointEll}.

\subsection{The second-order frame bundle}\label{sec2NdFrBund}

For any $L\in\LL(n,2n)$ we define the space\footnote{See, e.g., \cite[Chapter IV]{MR1202431}, \cite[Example 5.2]{kobayashi1996foundations}, for more details
on frame bundles.} of \emph{second-order frames} at $L$ as
\begin{equation}
\FF_L^2:=\{ j^2_0(\phi)\mid \phi:S^2V^*\rightarrow\LL(n,2n),\ \phi(0)=L,\ \phi\textrm{ local diffeomorphism at }0\},
\end{equation}
and the \emph{second-order frame bundle} of $\LL(n,2n)$ by
\begin{equation}
\FF^2:=\coprod_{L\in\LL(n,2n)}\FF^2_L.
\end{equation}
This bundle allows us to ``see'' the action of the positive-degree part of the group $P$. Recall that, by its definition \eqref{eqPrimaDefPi}, $P$ consists of
diffeomorphisms of $\LL(n,2n)$ preserving $\arctan(0)$. In particular, each $q\in P$ can be regarded as a local diffeomorphism $q:S^2V^*\rightarrow\LL(n,2n)$
sending 0 into $\arctan(0)$. Therefore, $j^2_0(q)\in\FF^2_{\arctan(0)}$ and we found the map
\begin{align}
P &\longrightarrow  \FF^2_{\arctan(0)},\nonumber\\
q&\longmapsto  j^2_0(q).\label{eqMapCheDovrebbeEssereIniettiva}
\end{align}
Obviously, $\FF^2_{\arctan(0)}$ is a group, and it can be proved that \eqref{eqMapCheDovrebbeEssereIniettiva} above is a group embedding. Therefore,
the structure group $P$ of the bundle $\Sp_{2n}\rightarrow\LL(n,2n)$ embeds into the structure group of the bundle $\FF^2\rightarrow\LL(n,2n)$.
Then the $P$-principal bundle $\Sp_{2n}\rightarrow\LL(n,2n)$ can be regarded as a \emph{reduction} of the second-order frame bundle of $\LL(n,2n)$.

The reduction is easier grasped on the Lie algebra level. Indeed, the Lie algebra of the group $\FF^2_{\arctan(0)}$ is
\begin{equation}
\gl(S^2V^*)\oplus (S^2(S^2V^*)\otimes S^2 V)
\end{equation}
and it contains $\gp$ as a subalgebra. The embedding is indicated by \eqref{eqAzioneInfinitesimaleMPuntoZero}. The 0-degree component of $\gp$ embeds naturally
into $\gl(S^2V^*)$. An element $B\in S^2V$, that is the 1-degree component of $\gp$, is mapped into the bilinear map
\begin{gather*}
S^2V^*\times S^2V^* \longrightarrow S^2 V, \\
(h,k)\longmapsto h\cdot B\cdot k.
\end{gather*}

Regarding $\Sp_{2n}$ as a sub-bundle of the second-order frame bundle of $\LL(n,2n)$ is an indispensable step when it comes to the problem of equivalence
of hypersurfaces in $\LL(n,2n)$. Such a problem is usually dealt with, in the spirit of Cartan, via  the moving frame methods, i.e., restrictions of $\FF^k$
to the embedded hypersurfaces at hand.

\subsection{Representation theory of $\Sp_{2n}$ and its subgroup $\GL_n$}\label{SecRepStuff}
The standard choice of a \emph{Cartan subalgebra} of $\sp_{2n}$ is given by the $n$-dimensional Abelian subalgebra
\begin{equation}
\gh:=\Span{\epsilon^i\otimes e_i\mid i=1,\dots, n}
\end{equation}
of diagonal matrices in $\gl(V)$. The \emph{fundamental weights} are then
\begin{equation}
\lambda_j:=\sum_{i=1}^je_i\otimes\epsilon^i ,\quad j=1,\dots,n
\end{equation}
(see \cite[\S 2.2.13]{MR2532439}) where $e_i\otimes\epsilon^i\in \gh^*$ is the basis element dual to $\epsilon^i\otimes e_i\in \gh$. For any $j=1,\dots,n$,
the fundamental weight $\lambda_j$ appears as the weight of the highest weight vector
\begin{equation}
v_{\lambda_j}:=e_1\wedge\ldots\wedge e_j
\end{equation}
in $\wedge^j(V\oplus V^*)$. Observe that
\begin{equation}
[v_{\lambda_n}]=\vol(\arctan(0)),
\end{equation}
that is, the Pl\"ucker image of the origin $\arctan(0)\in\LL(n,2n)$ is the line through the highest weight vector in $\wedge^n(V\oplus V^*)$. The subtle point
is that $\wedge^n(V\oplus V^*)$ is \emph{not} the highest weight module $V_{\lambda_n}$ of $v_{\lambda_n}$. Indeed, $\wedge^n(V\oplus V^*)$ is not irreducible
and
\begin{equation}\label{eqPluckEmbSpacRT}
V_{\lambda_n}=\wedge_0^n(V\oplus V^*)
\end{equation}
is precisely the subspace introduced in \eqref{eqSpazioMinimalePlucker} above. Therefore, $\p V_{\lambda_n}$ is the representation-theoretic way of describing
the Pl\"ucker embedding space.

Since $\GL(V)\subset\Sp(V\oplus V^*)$, any irreducible representation of $\Sp(V\oplus V^*)$ becomes a (not necessarily irreducible) representation of $\GL(V)$.
In particular, the irreducible $\Sp(V\oplus V^*)$-representation $V_\lambda$ with highest weight $\lambda$ splits into several $\GL(V)$-irreducible
representations. Only one of the latter contains the weight vector $v_\lambda$ and therefore it will be denoted by $L_\lambda$:
\begin{equation*}
L_\lambda=\textrm{the unique }\GL(V)\textrm{-irreducible component of }V_\lambda\textrm{ containing } v_\lambda.
\end{equation*}
For instance, $V_{\lambda_1}$ is the $2n$-dimensional fundamental representation $V\oplus V^*$, whereas $L_{\lambda_1}$ is simply $V$. Similarly,
$V_{\lambda_n}$ is the (very large) de-projectivised Pl\"ucker embedding space for $\LL(2n,n)$, whereas  $L_{\lambda_n}$ is the one-dimensional line
$\vol(\arctan(0))$!

\subsection{The tautological line bundle and $r$-th degree hypersurface sections}\label{SecRThDegreeSection}

Having at one's disposal the $P$-principal bundle \eqref{eqPrincBundleSpn} and regarding the $\GL_n$-irreducible representation $L_\lambda$ as a representation
of $P$, one can form the associated vector bundle
\begin{equation}
\Sp_{2n}\times^PL_\lambda.
\end{equation}
For instance, with $\lambda=\lambda_1$ one obtains the tautological bundle introduced earlier (cf. \eqref{eqTautBundPRE},  \eqref{eqTautBund}) and with
$\lambda=\lambda_n$ one obtain the tautological \emph{line} bundle $\OO_{\LL(n,2n)}(-1)$. The readear should bear in mind that the former is a rank-$n$ bundle,
whereas the latter has rank 1. Indeed, there are two ``tautological'' principles at play here, following from the fact that $\LL(n,2n)$ is made of
$n$-dimensional linear subspaces, and $\p V_{\lambda_n}$ consists of lines, respectively.

By definition, $\OO_{\LL(n,2n)}(-1)$ is the pull-back via the Pl\"ucker embedding of the   tautological line bundle over the Pl\"ucker embedding space
$\p V_{\lambda_n}$. Indeed, the fibre of $\OO_{\LL(n,2n)}(-1)$ at $L\in\LL(n,2n)$ is $\vol(L)$ itself, understood not as a point of $\p V_{\lambda_n}$,
but as an abstract one-dimensional linear space. Such is the standard notation of Algebraic Geometry: over the projective space $\p(W)$ one has a group
(isomorphic to $\Z$) of linear bundles
\begin{equation*}
\OO(-1)_{[w]}:=\Span{w},\quad \OO(1) :=\OO(-1) ^*,\quad  \OO(\pm r) :=\OO(\pm 1)^{\otimes r},\quad \forall [w]\in\p W,\,r\in\N.
\end{equation*}
Given a hyperplane $\pi\in\p V_{\lambda_n}^*$, we may form the \emph{hyperplane section} $\Sigma_{\pi}:=\pi\cap\LL(n,2n)$. The same hypersurface
$\Sigma_{\pi}$ can be described as the zero locus of a suitable section of $\OO_{\LL(n,2n)}(1)$, the dual of $\OO_{\LL(n,2n)}(-1)$. Indeed,
let $\alpha\in V_{\lambda_n}^*$  be a linear form such that $\pi=\p(\ker\alpha)$. Then $\alpha$ can be restricted to each line $\vol(L)$, thus yielding
a section (still denoted by $\alpha$) of $\OO_{\LL(n,2n)}(1)$. The value of the section $\alpha$ at $L$ is simply $\alpha|_{\vol(L)}$. Therefore, the zero
locus of $\alpha$ is made precisely by those $L$, such that $\vol(L)\subset\ker\alpha$, that is, $\vol(L)\in\pi$, which is precisely $\Sigma_{\pi}$.

A central question in the geometry of PDEs is:  \emph{how to recognise a hyperplane section?} In the above language of induced bundles, this is the same
as asking: \emph{when a section of $\OO_{\LL(n,2n)}(1)$ comes from a linear form $\alpha\in V_{\lambda_n}^*$?}

In general, the map
\begin{equation}\label{eqPrimaFrecciaBGG}
S^rV_{\lambda_n}^*\longrightarrow \Gamma\bigl(\LL(n,2n),\OO_{\LL(n,2n)}(r)\bigr)
\end{equation}
associating with a degree-$r$ homogeneous polynomial on $V_{\lambda_n}$ a (global) section of  $\OO_{\LL(n,2n)}(r)$, the $r$-th power of
$\OO_{\LL(n,2n)}(1)$, can be \emph{resolved}. More precisely, there exists a differential operator $\square^r$ acting on sections of $\OO_{\LL(n,2n)}(r)$,
whose kernel is precisely the image of \eqref{eqPrimaFrecciaBGG}. The construction of  $\square^r$ is by no means trivial and it is based on the so-called
BGG resolution \cite[Theorem 5.9]{MR3603758}.

From the point of view of PDEs, the operator $\square^r$ is to be understood as a \emph{test}, that is, as a criterium to establish whether a given
second-order PDE $F(p_{ij})=0$ belongs to the well-defined ($\Sp_{2n}$-invariant) class of $r$-th degree hypersurface sections. Running the test entails
applying $\square^r$ to the function $F$ defining the equation at hand. Therefore, the above class of second-order PDEs is to be understood as the set of
solutions of the special differential equation $\square^r(F)=0$. The same   idea will be applied  to another important class of PDEs, the integrable ones,
see Section  \ref{SecInt} below.

\subsection{Rank-one vectors}\label{secRO}
An immediate consequence of the fundamental isomorphism \eqref{eqIdentFund} is that the projectivized tangent bundle $\p T\LL(n,2n)$ contains a proper
sub-bundle, namely the bundle
\begin{equation}
\sR:=\{[v]\in  \p T\LL(n,2n)\mid \rank v=1\}
\end{equation}
of (projective classes of) \emph{rank-one vectors}. Indeed, \eqref{eqIdentFund} allows us to speak of the   rank of the vector $v$, and such a notion is
well-defined and depends only on the projective class of $v$.

We provide now an interesting characterisation of rank-one tangent vectors. Let $L\in\LL(n,2n)$ be an arbitrary point, and $v$ a tangent vector at $L$
to $\LL(n,2n)$. Let $\gamma=\gamma(\epsilon)$ be a curve passing through $L$ with speed $v$. Then, each point $\gamma(\epsilon)$ can be interpreted as
a Lagrangian subspace of  $V\oplus V^*$, and in particular $\gamma(0)=L$. Observe that, even for  small values of $\epsilon$, the intersection
$\gamma(\epsilon)\cap\gamma(0)$ needs not to be nontrivial.

Here it comes the peculiarity of rank-one vectors: $v$ is rank-one if and only if the curve $\gamma$ can be chosen in such a way as the intersection
$\gamma(0)\cap\gamma(\epsilon)$ is a \emph{fixed hyperplane} $\Sigma\subset L$ (that is, not depending on $\epsilon$) for all $\epsilon$ in a small
neighbourhood of zero. In other words, there is a correspondence between rank-one tangent vectors at $L$ and hyperplanes $\Sigma\subset L$, that is,
elements of $\p L^*$. In one direction, such a correspondence is quite evident.

\abs-0.1 \aabd-1
Thanks to the homogeneity of $\LL(n,2n)$  we can work at the origin $\arctan(0)$ (see Section \ref{secHomStructLG} above).
Let $[\alpha]\in\p \arctan(0)^*$ represent the hyperplane $\Sigma_\alpha:=\ker\alpha$.
We need to describe the generic Lagrangian plane $L$, which is ``close'' to  $\arctan(0)$, and intersects  the latter precisely along $\Sigma_\alpha$.
Since $L$ has to be ``close'' to $\arctan(0)$, we can assume it to belong to the big cell, that is, to be of the form $\arctan(h)$, for some $h\in S^2V^*$.
The key remark, rather obvious, is that
\begin{equation}
\arctan(0)\cap\arctan(h)=\ker h.
\end{equation}
So, the above intersection coincides with $\Sigma_\alpha$ if and only if $\ker\alpha=\ker h$. That is, if and only if the quadratic form $h$ is proportional
to the square of the linear form $\alpha$. Therefore, the correspondence between hyperplanes in $L$ and rank-one tangent vectors to $\LL(n,2n)$ at $L$ is
nothing but
\begin{align}
\p L^* &\longrightarrow \p S^2L^*\, \nonumber\\
\vspace{3pt}[\alpha] &\longmapsto [\alpha^2],\label{eqSegre}
\end{align}
one of the most fundamental maps in classical Algebraic Geometry: the Veronese embedding \cite[Example 2.4]{harris1992algebraic}. The above map
\eqref{eqSegre},  in the context of second order PDEs, allows us to establish an important relationship between objects depending on second order derivatives
(elements of $S^2L^*$ are reminiscent of Hessian matrices) and objects depending on first order derivatives (elements of $L^*$ correspond to covectors on the
space of independent variables). This point of view will be clarified in Section \ref{secGoursatMA}.

\subsection{Characteristics}\label{SecChar}
In compliance with the terminology found, e.g., in  \cite[Formula (3.1)]{MR2985508},  \cite[Chapter VI]{MR1083148} and \cite[Part II]{Smith2019}, we denote by
\begin{equation}\label{eqPrologSigma}
\Sigma^{(1)}:=\{L\in\LL(n,2n)\mid L\supset\Sigma\}
\end{equation}
the \emph{prolongation} of the hyperplane $\Sigma\in\p L^*$. As we have already pointed out, $\Sigma^{(1)}$ is a line passing through $L$ itself (see
Section \ref{secRO}). In fact, via Pl\"ucker embedding, $\Sigma^{(1)}$ becomes an actual projective line in $\p V_{\lambda_n}$, see, e.g.,
\cite[Proposition 2.1]{MR2985508}.

Moreover, if $\Sigma=\ker\alpha$, with $\alpha\in L^*$, then
\begin{equation}
T_L\Sigma^{(1)}=\Span{\alpha^2},
\end{equation}
as a subset of $S^2L^*=T_L\LL(n,2n)$ (recall formula \eqref{eqIdentFundPointEll}). Let $\E\subseteq \LL(n,2n)$ be a submanifold and $L\in\E$. Then $\Sigma$
is called a \emph{characteristic} (resp., \emph{strong characteristic})  for $\E$ at $L$ if $\Sigma^{(1)}$ is tangent to (resp., contained in) $\E$.
These notions will be essential in the analysis of the well-posedness of initial value problems for PDEs, see Section \ref{secCharVar}.

\subsection{The Lagrangian Chow transform}\label{secLagChowTr}
So far we have worked with Lagrangian---i.e., \emph{maximally isotropic}---subspaces $L$ of $V\oplus V^*$. The hyperplanes $\Sigma$ appearing in
Section \ref{SecChar} are the first instances of \emph{sub-maximal} isotropic subspaces (in this case, ($n-1$)-dimensional). In fact, nothing forbids
considering the sets
\begin{equation}
\LL(i,V\oplus V^*):=\{L\in\Gr(i,V\oplus V^*)\mid \omega|_L\equiv 0\},\quad i=1,\dots, n,
\end{equation}
which we may call ``Lagrangianoid Grassmannians'', as well as the corresponding \emph{incidence correspondences}:
\begin{equation}\label{eqDefFlagCheArmandoTantoDesiderava}
\Fl^\iso(V\oplus V^*; i,j):=\{ (L_1,L_2)\in \LL(i, V\oplus V^*)\times\LL(j, V\oplus V^*)\mid L_1\subset L_2\},
\end{equation}
where $1\leq i<j\leq n$.

\abs0.2 \aabd2
In particular, an important  role is played by $\LL(1,V\oplus V^*)\equiv\p (V\oplus V^*)$ and $\LL({n-1},V\oplus V^*)$. Indeed, in the geometric theory of PDEs,
the former describes rank-one subdistributions of  the contact distribution, and the latter describes infinitesimal Cauchy data. The two notions coincide
for $n=2$.\footnote{The classical reference in the book \cite{Petrovsky1992}. Different treatments of the subject, sometimes closer in  spirit to the present
paper, can be found, e.g., in \cite{Vitagliano2013,MorenoCauchy,MR2985508}.}

A classical observation in Algebraic Geometry is that all these ``Lagrangianoid Grassmannians'' are tied together by means of the  incidence correspondences
\eqref{eqDefFlagCheArmandoTantoDesiderava}. Indeed, the above-defined sets of \emph{isotropic} flags  fit into the following double fibration:
\begin{equation}\label{eqDiagDoubFibr1}
\gathered
\xymatrix{
 &\Fl^\iso(V\oplus V^*; i,j)\ar[dl]_{p_i}\ar[dr]^{p_j}&\\
 \LL(i, V\oplus V^*) && \LL(j, V\oplus V^*),
}
\endgathered
\end{equation}
with $i<j$. For instance, with $i=n-1$ and $j=n$   diagram \eqref{eqDiagDoubFibr1} reads
\begin{equation}\label{eqDiagDoubFibrNMenoUno}
\gathered
\xymatrix{
 &\Fl^\iso(n-1,n, V\oplus V^*)\ar[dl]_{p_{n-1}}\ar[dr]^{p_n}&\\
 \LL(n-1,V\oplus V^*) && \LL(n,V\oplus V^*),
}
\endgathered
\end{equation}
and for any $\Sigma\in \LL(n-1,V\oplus V^*) $, the ``double fibration transform'' $p_n(p_{n-1}^{-1}(\Sigma))$ of $\Sigma$ is precisely the prolongation
$\Sigma^{(1)}$ defined by \eqref{eqPrologSigma}. Conversely, for any $L\in \LL(n,V\oplus V^*) $, the ``inverse double fibration transform''
$p_{n-1}(p_{n}^{-1}(L))$ of $L$ is nothing but $\p L^*$.

Another interesting example is obtained with $i=1$ and $j=n$. Diagram  \eqref{eqDiagDoubFibr1} then reads
\begin{equation}\label{eqDiagDoubFibrUno}
\gathered
\xymatrix{
 &\Fl^\iso(V\oplus V^*;1,n)\ar[dl]_{p_{1}}\ar[dr]^{p_n}&\\
 \p(V\oplus V^*) && \LL(n,V\oplus V^*).
}
\endgathered
\end{equation}
The above diagram allows us to recast a simple but useful theorem, known in Algebraic Geometry as the \emph{Chow form/transform}: if $X\subset
\p(V\oplus V^*)$ is a smooth variety of dimension $n-1$, then its ``double fibration transform'' $p_n(p_{1}^{-1}(X))$ is a smooth hypersurface in
$\LL(n,V\oplus V^*)$, of the \emph{same degree} as $X$ \cite[Lemma 23]{AGMM_CommContMath}. The latter will be referred to as the
\emph{Lagrangian Chow transform} of $X$. We stress that the notion of \emph{degree} in  $\LL(n,V\oplus V^*)$ refers to the surrounding Pl\"ucker embedding
space.

As a  nice example consider  an $n$-dimensional (not necessarily Lagrangian) subspace $D\subset V\oplus V^*$. Then $\p D$ is a (smooth) $(n-1)$-dimensional
variety in $ \p(V\oplus V^*)$ whose Lagrangian Chow transform reads
\begin{equation}\label{eqProtoGoursatMAE}
\det(D-h)=0,
\end{equation}
where $D$ is the (not necessarily symmetric) $n\times n$ matrix corresponding to the subspace $D$ in the big cell $V^*\otimes V^*$ of $\Gr(n,V\oplus V^*)$
and $h$ is the symmetric  $n\times n$ matrix corresponding to the generic element of the big cell $S^2V^*$ of $\LL(n,V\oplus V^*)$. Observe that
\eqref{eqProtoGoursatMAE}, though containing all the minors of $h$, is \emph{linear} in the Pl\"ucker coordinates, as predicted by the theorem.

\abs-0.1 \aabd-1
The second order PDEs corresponding to the Lagrangian Chow transforms of the $n$-di\-men\-sion\-al sub-distributions of the contact distribution are the
so-called \emph{Goursat-type Monge--Amp\`ere equations}, introduced by E. Goursat in 1899  \cite{MR1504329}, way before the inception of the Chow transform,
see Section \ref{secGoursatMA} below. It is precisely thanks to the introduction of the Lagrangian Chow transform that the notion of a Goursat-type
Monge--Amp\`ere equation can   be generalised to arbitrary  \emph{conic sub-distributions} of the contact distribution \cite[Section 9]{Buczyski2019}.

\subsection{A few remarks on $\LL(2,4)$ and $\LL(3,6)$}
We conclude this  survey of the rich geometry of $\LL(n,2n)$ by pointing out the peculiarities of two low-dimensional examples, namely when $n=2$ or $n=3$.
The case $n=2$ is examined from top to bottom in the paper \cite{MR2876965}. Even if the case $n=3$  does not  boast its own treatise, the reader will find
specific facts and results in \cite[\S 5]{Russo2019} and \cite[\S 4.2]{MR3603758}. We do not review here all that can be found in the aforementioned works---we
rather highlight the  origin of the interestingness and diversity of these two cases.

The departing point is the fact, already pointed out, that $\LL(n,2n)$ is always  isotropic with respect to the natural two-form defined on the
(de-projectivised) Pl\"ucker embedding space $V_{\lambda_n}$, see \eqref{eqPluckEmbSpacRT}.

In the case $n=2$, this two-form is symmetric (see \eqref{eqFormSimmPluckSpace}) and we denote it here by~$g$. Therefore, since the codimension of $\LL(2,4)$
in  $\p V_{\lambda_2}\equiv\p^4$ is one, $\LL(2,4)$ must coincide with the (projectivised) null cone of $g$ in $\p^4$.

In the real case, $g$ has signature $(+++--)$ and then $\LL(2,4)$ inherits a conformal structure of signature $(++-)$. Such an  $\Sp_4(\R)$-invariant conformal
structure is precisely the one that has been used by The to carry out a classification of hypersurfaces in $\LL(2,4)$ by the method of moving frames
\cite{MR2876965}. The same structure has also been used by the authors to characterise the hyperplane sections of $\LL(2,4)$ in terms of the trace-free second
fundamental form \cite[Corollary 4.2]{MR3603758}.

Another peculiarity of $\LL(2,4)$ which is worth recalling is that  $\LL(2,4)$ is isomorphic to the so-called \emph{Lie quadric}. This is the moduli space
of \emph{all} circles in $\R^2$, i.e., including also those with zero or infinite radius. Such an isomorphism was essentially known to S.~Lie himself
\cite[Lie's Memoir on a Class of Geometric Transformations, \S 9]{MR0106139}, though it can be rephrased in modern language by using Hopf fibration,
see  \cite[\S 5]{2010arXiv1002.0243A} and~\cite{MO165717}.

Passing to the case $n=3$, we see that the 6-dimensional $\LL(3,6)$ does not carry any natural conformal structure in the usual sense. Nevertheless a
``trivalent'' analogue of a conformal structure can still be defined on $\LL(3,6)$. Such a structure has been exploited by the authors to characterise
hyperplane sections of $\LL(3,6)$ in terms of a suitable generalisation of the  trace-free second fundamental used in the case $n=2$ \cite[Section~4.2]{MR3603758}.

Another really intriguing feature of $\LL(3,6)$, or rather of its Pl\"ucker embedding in~$\p^{13}$, is that such an embedding can be regarded as an
appropriate generalisation of the twisted cubic in $\p^3$, whereby the field of complex number has been replaced by the Jordan algebra of symmetric $3\times 3$
matrix. This analogy played a fundamental role in a recent analysis of PDEs with prescribed group of symmetries \cite{The2018}. A gentle introduction to it can
be found in \cite[\S 5]{Russo2019}.

\section{Hypersurfaces in the (real) Lagrangian Grassmannian  and second order PDEs}\label{secPart2}
In the second part of this paper we examine more in depth the geometry of hypersurfaces in the Lagrangian Grassmannian $\LL(n,2n)$. Some of the key notions,
like those of a \emph{hyperplane section} (Section \ref{secDualVAr}), of an \emph{$r$-th degree section} (Section \ref{SecRThDegreeSection}) and of the
\emph{characteristic of a hypersurface} (Section \ref{SecChar}), have already been introduced above. It was also anticipated that these ideas were going
to have interesting incarnations in the context of second order PDEs. All of this will be explained below.

From now on, we work in the real smooth category.

\subsection{Contact manifolds and second order PDEs}
The idea of framing second order PDEs against the general background of contact manifolds and their prolongations is rather old and, in a sense, it belongs
to the mathematical folklore. An excellent treatise of this topic is the book \cite{MR2352610} though a slenderer introduction can be found
in \cite{EMMS_CompMan2018}.

The departing point is a contact manifold $(M,\CC)$, that is a $(2n+1)$-dimensional smooth manifold equipped with a one-codimensional distribution $\CC$,
such that the Levi form,
\begin{gather}
\omega:\CC\wedge\CC \longrightarrow  TM/\CC,\nonumber\\
(X,Y) \longmapsto [X,Y]+\CC,
\end{gather}
is non-degenerate. The so-called \emph{Darboux coordinate} may help to  clarify the picture: $M$~is (locally) described by the coordinates
\begin{equation}\label{eqDarbCoord}
   (x^1,\dots, x^n,u,p_1,\dots, p_n),
\end{equation}
the distribution $\CC$ is (locally) spanned by the $2n$ vector fields
\begin{equation*}
\underbrace{D_1:=\partial_{x^1}+p_1\partial_u,\dots, D_n:=\partial_{x^n}+p_n\partial_u}_{\textrm{``total'' derivatives}},
\underbrace{\partial_{p_1},\,\dots,\partial_{p_n}}_{\textrm{``vertical'' derivatives}},
\end{equation*}
and  (locally)
\begin{equation}
 \omega=dx^i\wedge dp_i\, .
\end{equation}
The next step consists in regarding each contact plane $\CC_p$, with $p\in M$, as a symplectic linear space (thanks to the symplectic form $\omega_p$)
and in constructing the corresponding Lagrangian Grassmannian $\LL(\CC_p):=\LL(n,\CC_p)$.   One readily verifies that the total derivatives and the vertical
derivatives are dual to each other via $\omega$, that is, they can be identified with the vectors $e_i$ and the covectors $\epsilon^i$ introduced
in Section~\ref{secPrelim}, respectively. Then, following the same  procedure as in Section~\ref{secDefLagGrass},  we obtain coordinates $p_{ij}$ on
$\LL(\CC_p)$. Doing the same for any point $p$ one obtains a bundle
\begin{equation}
\LL(\CC):=\coprod_{p\in M} \LL(\CC_p)\longrightarrow M
\end{equation}
with fibre coordinates $p_{ij}$, known as the \emph{Lagrangian Grassmannian bundle} of $M$ or the \emph{first prolongation} of $M$ and sometimes denoted
by $M^{(1)}$.

Then a hypersurface $\E\subset M^{(1)}$, being locally represented as
\begin{equation}\label{eqPDE}
 \E: F(x^1,\dots, x^n,u,p_1,\dots, p_n,\dots, p_{ij},\dots)=0,
\end{equation}
clearly corresponds to a \emph{second order   PDE}. Perhaps it is less evident that a \emph{solution} of $\E$ is captured by a Lagrangian submanifold
$L\subset M$, such that its tangent lift $TL$ is contained into $\E$. In the  coordinates \eqref{eqDarbCoord} of Darboux, $L=L_f:=\{\x, f(\x), \dots,
(\partial_{x^i}f)(\x),\dots \}$, where $f$ is a function of $\x=(x^1,\dots, x^n)$, and it is not hard to prove that $TL$ (the set of all the tangent
$n$-dimensional subspaces to $L$) coincides with
\begin{equation}\label{eqElleEffeUno}
 L_f^{(1)}:= \{\x, f(\x), \dots, (\partial_{x^i}f)(\x),\dots, \dots, (\partial_{x^i}\partial_{x^j}f)(\x),\dots \}\subset M^{(1)},
\end{equation}
so that $ L_f^{(1)}\subset\E$ if and only if the function $f$ fulfills the (familiar looking) PDE appearing in \eqref{eqPDE}.

From now on we make the (non-restrictive) assumption that $\E$ is actually a sub-bundle of $M^{(1)}$. Then the fibres $\E_p$ of $\E$ are hypersurfaces
in the corresponding Lagrangian Grassmannians $\LL(\CC_p)$, with $p\in M$. So,  we are in position of utilising the theoretical machinery developed in the
first part. Essentially, we are going to work with a \emph{family} of symplectic spaces, Lagrangian Grassmannians and hypersurfaces of the latter, rather
than with a fixed one. Besides the appearance of a fancy index ``$p$'', the techniques remain unchanged.

A subtler point, which may have escaped the hasty reader, is that passing from the point-wise perspective (``microlocal'', as some love to say) to the global
framework, the equivalence group has changed from the finite-dimensional Lie group $\Sp_{2n}$ to the infinite-dimensional \emph{contact group} $\Cont(M)$.

\subsection{Nondegenerate second order  PDEs and their symbols}\label{secNonDeg}
If one's ultimate goal is to be able to   setup the equivalence problem for second order  PDEs, then there is one rough distinction that can be made from the
very beginning.

A hypersurface $\E\subset\LL(n,2n)$ is called \emph{non-degenerate at $L$} if the tangent hyperplane $T_L\E$, understood as a line in $S^2L$ via the dual of
identification \eqref{eqIdentFund} is made of non-degenerate elements. Then $\E$ is called \emph{non-degenerate} if it is non-degenerate at all points.
Finally, a second order  PDE $\E\subset M^{(1)}$ is \emph{non-degenerate} if so  are all its fibres. Obviously, the property of being non-degenerate is
$\Cont(M)$-invariant and hence defines a well-behaved class of second order PDEs.

The fundamental correspondence \eqref{eqIdentFund} reads now, in terms of the  local Darboux coordinates \eqref{eqDarbCoord},
\begin{gather}
S^2 L^*  \longrightarrow  T_L\LL(n,2n),\nonumber\\
\dd x^i\odot \dd x^j \longmapsto \partial_{p_{ij}}|_L\, .
\end{gather}
Therefore, if $\E=\{F=0\}$ is a hypersurface in  $\LL(n,2n)$, then $\dd_LF$ can be regarded as an element of $S^2L$, viz.
\begin{equation}\label{eqSymbolF}
\dd_LF=\frac{\partial F}{\partial p_{ij}}\biggr|_L\dd p_{ij}\longleftrightarrow \frac{\partial F}{\partial p_{ij}}\biggr|_L \partial_{x^i}\odot \partial_{x^j}\,.
\end{equation}
The symmetric rank-two contravariant tensor appearing at the right-hand side of \eqref{eqSymbolF} is of paramount importance in the theory of PDEs. It is
called the \emph{symbol} of $F$ at $L$. If the dependence upon $L$ is discarded then one has a section of the bundle $S^2L|_\E$ (beware of the syncretism
of the symbol $L$, cf. \eqref{eqIdentFund} and \eqref{eqIdentFundPointEll}), still called the \emph{symbol} of $F$. Finally, if $\E\subset M^{(1)}$ is a second
order  PDE, then the \emph{symbol} of $F$ must be understood as a section of a bundle over $\E$, whose restriction to the fibre $\E_p$ is the aforementioned
bundle $S^2L|_{\E_p}$. Such a proliferation of ``bundles upon bundles'' is a congenital feat of the theory and the reader must cope with it, see also
\cite[Section 3]{Smith2019}. Using the same symbol for the various incarnations of the same concept, far from bringing in more confusion, is the only way
to keep the notation bearable.

Now we must face a fundamental problem in the theory of hypersurfaces, that is the fact that $F$ is not, of course, uniquely determined by $\E$ and it is
$\E$ that we wish to study, not $F$. Usually things are simpler with $F$, but then one has to ensure the result to be independent upon the
choice\footnote{Borrowing a terminology from Algebraic Geometry, we call the \emph{ideal} of $\E$ the ideal in $C^\infty(M^{(1)})$ of functions vanishing on
$\E$.} of $F$ in the ideal of $\E$. Another way out is to prove results directly on $\E$, but this usually demands a deeper abstraction.

For instance, the \emph{symbol of the equation $\E$} at $L\in M^{(1)}$ is the element
\begin{equation}\label{eqDefSimbProject}
 \Smbl_L(\E):=[\dd_L F]\in \p (S^2 L),
\end{equation}
whereas the previously defined symbol of $F$ at $L$ is just a representative of it. One is more conceptual, the other more treatable. Nevertheless, both allow
us to rephrase the notion of non-degeneracy: the PDE $\E$ is non-degenerate at $L$ if its symbol at $L$ is a generic element of  $\p (S^2 L)$ or, equivalently,
if the symbol of \emph{any representative} $F$ of $\E$ is a non-degenerate rank-two symmetric tensor on $L$.

\subsection{Symplectic second order PDEs}

The various versions of the above notion of non-degeneracy (in a point, in a fibre, everywhere) stressed the main issue of passing from the study of
hypersurfaces in $\LL(n,2n)$ to the study of second order PDEs $\E\subset M^{(1)}$:   the fibres of $\E$ may fulfill some special property (e.g., that of
being non-degenerate) over some subset $U\subset M$ and, simultaneously,  may \emph{not} fulfill it over $M\setminus U$. This is the main source of
additional difficulties: two equations of ``mixed type'' may not be $\Cont(M)$-equivalent for topological reasons (e.g., because the locus $U$ of the first
equation is not homeomorphic to the analogous locus of the second equation).

A reasonable compromise between the $\Sp_{2n}$-equivalence problem and the $\Cont(M)$-equivalence problem is provided by the sub-class of second order PDEs
that locally look like
\begin{equation}\label{eqPDE_Hirota}
 \E: F(\dots, p_{ij},\dots)=0,
\end{equation}
that is, exactly like \eqref{eqPDE}, but without explicit dependency upon $x^1,\dots, x^n,u,p_1,\dots, p_n$.
Such a class will be called the class of
\emph{symplectic second order PDEs} in compliance with  the terminology adopted, e.g., in \cite{DoubrovFerapontov,MR2765729, MR2805306}.\footnote{According
to another school, this is the class of \emph{Hirota-type} second order PDEs, see e.g., \cite{doi:10.1093/imrn/rnp134,2017arXiv170708070F}.}  More
geometrically, one can speak about symplectic second order PDEs when the bundle $M^{(1)}\longrightarrow M$ is trivial, i.e.,  $M^{(1)}=M\times\LL(n,2n)$
and $\E$ is the pull-back of a hypersurface (still denoted by $\E$) in $\LL(n,2n)$. Hence, the modifier ``symplectic'' alludes to the fact that the equivalence
group is still  $\Sp_{2n}$, even though the equation is defined over the contact manifold $M$. From now on, unless otherwise stated, all second order PDEs
are assumed to be (everywhere) non-degenerate and symplectic. The same symbol $\E$ will be used both for the sub-bundle of~$M^{(1)}$ and for an its generic
fibre. The context will help the reader to know which is which.

It may happen that very hard questions for general   second order PDEs become almost trivial in the context of symplectic second order PDEs. For instance,
the problem of linearisability of a general parabolic Monge--Amp\`ere equation, up to contactomorphisms, was raised by R. Bryant \cite{MR1334205} and to date
it is still open, whereas its analogue for symplectic Monge--Amp\`ere equations is (relatively) trivial, see \cite[Theorem 1.4]{Russo2019}. Obviously, the class
of symplectic second order PDEs is \emph{not} $\Cont(M)$-invariant.

\subsection{The characteristic variety}\label{secCharVar}

Before introducing the simplest yet nontrivial class of PDEs, we recast the notion of a characteristic  in the present context of PDEs. Recall that,
for any point $L\in\E$, a hyperplane $\Sigma\in\p L^*$ is called a \emph{characteristic} (resp., \emph{strong characteristic}) for $\E$ at $L$ if the rank-one
line $\Sigma^{(1)}\subset T_L\LL(n,2n)$ it tangent to $\E$ at $L$ (resp., contained into $\E$), see Section \ref{SecChar}. Let now $\E\subset M^{(1)}$ be a PDE,
and $L\in\E$. The subset
\begin{equation}\label{eqDefCharVar}
  \Xi_L(\E):=\{ \Sigma\in \p L^*\mid \Sigma^{(1)}\textrm{ is tangent to }\E\textrm{ at }L\}\subset \p L^*
\end{equation}
is called the \emph{characteristic variety} of $\E$ at $L$. Their (disjoint) union, for all $L\in\E$, forms a bundle over $\E$ called simply the
\emph{characteristic variety} and denoted by $ \Xi(\E)$.

The conceptual definition \eqref{eqDefCharVar} may be abstruse, but Darboux coordinates make it easily accessible to computations. It is easy to see that
\eqref{eqDefCharVar} can be equivalently formulated as
\begin{equation}\label{eqDefCharVar2}
 \Xi_L(\E)=\{[\alpha]\mid \alpha\in L^*,\ \Smbl_L(F)(\alpha^2)=0\}.
\end{equation}
Here  $\Smbl_L(F)$ is the symbol of $F$ at $L$, as in the right-hand side of \eqref{eqSymbolF}. Observe that the condition at the right-hand side of
\eqref{eqDefCharVar2} is  independent upon the choice of $F$ in the ideal of $\E$.

In this section we merely provide the definition of the characteristic variety $ \Xi(\E)$. A~careful examination of all the properties of $ \Xi(\E)$ and
ramifications would fill a separate treatise. For more information, we refer the reader to \cite{Smith2019} in this very volume and to \cite{Vitagliano2013}
and references therein. We just make  two final remarks.

First, the characteristic variety $ \Xi(\E)$ can be used to carry out a rough classification of PDEs. For instance, $\E$ is non-degenerate at $L$ iff
$\Xi_L(\E)$ is a non-degenerate\footnote{Beware that non-degenerate does not mean non-\emph{irreducible}.} quadric. Similarly, $\E$ is elliptic at $L$ iff
$\Xi_L(\E)$ is empty. In the case $n=2$, $\E$ is parabolic at $L$ iff $\Xi_L(\E)$ consists of two lines. And this list of examples may continue.

Second,  the characteristic variety $\Xi(\E)$  plays a fundamental role in the initial value problem. Assume, to make things even simpler, that a
characteristic $\Sigma$ is strong. Then the entire line $\Sigma^{(1)}$ is contained into $\E$. This means that there is a family, parametrised by $\p^1$,
of infinitesimal solutions to $\E$ \emph{admitting the same initial $($infinitesimal\/$)$ datum} $\Sigma$. In other words, if the initial datum is tangent to
$\Sigma$ (in which case the initial datum is called \emph{characteristic}), then the Cauchy--Kowalewskaya theorem fails in uniqueness. More examples
clarifying this property of  $ \Xi(\E)$   can be found in the above-cited paper \cite{Vitagliano2013}.

\subsection{Hyperplane sections and PDEs of Monge--Amp\`ere type}\label{secHypMAEs}
In the literature, \emph{the} Monge--Amp\`ere equation  is usually understood to be
\begin{equation}\label{eqTheMongeAmpere}
   \det(p_{ij})=f(x^i,u,p_i).
\end{equation}
It is at the very heart of a feverish research activity: for instance, the book \cite{MR1964483}, concerning the problem of the \emph{optimal mass
transportation},  gathered almost one thousand citations in a dozen of years. Besides countless scientific and technological applications, the problem of
optimal mass transportation can be formulated in important economical models, in the form of \emph{optimal allocation of resources}. This led, among many other
things, to a Nobel prize in the economic sciences for Kantorovich \cite{MR3220439,harvey1982,MR2945626,EMMS_CompMan2018}.

On a more speculative level,  one can ask for which functions $f$ in \eqref{eqTheMongeAmpere} one obtains a
($\Cont(M)$-invariant) class of PDEs. For instance, there exists a family of functions $f$  such that the corresponding equation \eqref{eqTheMongeAmpere}
can be brought  into the linear form\footnote{The simplest case of such a function is $f=0$.}
\begin{equation}
 p_{11}=0,
\end{equation}
by means of a (partial)
Legendre transformation (that is a particular element of  $\Cont(M)$), see, e.g., \cite{EMMS_CompMan2018,MR2985508}.
This is the easiest example of a $\Cont(M)$-invariant subclass of Monge--Amp\`ere equations---those having an \emph{integrable} characteristic distribution.
A~linear (symplectic) second order PDE
\begin{equation}\label{eqCretinaEqLineare}
 \E:F(p_{ij})= B^{ij}p_{ij}=0,\quad B^{ij}\in\R,
\end{equation}
is such that its representative $F$ fulfills the system of second-order PDEs
\begin{equation}\label{eqCretinaLinearita}
 \frac{\partial^2 F}{\partial p_{ij}\partial p_{hk}}=0,\quad \forall i,j,h,k.
\end{equation}
Equation \eqref{eqCretinaLinearita}, that is a PDE imposed on the left-hand side of another PDE (in this case, $\E$), is what we shall call a \emph{test}
later on. The key feature of \eqref{eqCretinaLinearita} is that it is \emph{not} $\Cont(M)$-invariant. Making \eqref{eqCretinaLinearita} into a
$\Cont(M)$-invariant test is not an easy task, and the heavy machinery used in \cite{MR3603758} confirms that; see also \cite{ Moreno2017}.
Nevertheless, the result is surprisingly simple, and even easy to guess. If we declare that $\E$ \emph{passes the Monge--Amp\`ere test} if and only if
\begin{equation}\label{eqMongAmpTestBrutal}
 \frac{\partial^2 F}{\partial p_{(ij}\partial p_{hk)}}=0,\quad \forall i,j,h,k,\quad\textrm{for some representative }F\textrm{ of }\E,
\end{equation}
then this test \emph{is} $\Cont(M)$-invariant. The curious reader may run it on \eqref{eqTheMongeAmpere} just for fun.

In the paper \cite{MR3603758} the authors have proved that a (symplectic) second order PDE $\E$ passes the Monge--Amp\`ere test if and only if $\E=\{F=0\}$,
with
\begin{equation}\label{eqMultDimMAE}
 F=A+B^{ij}p_{ij}+C^\bullet(2\times2\textrm{ minors})+\ldots+D^{ij}p_{ij}^\#+D\det(p_{ij}),
\end{equation}
which coincides with the classical definition of a (general) Monge--Amp\`ere equation with constant coefficients, see, e.g. \cite[Formula (0.5)]{MR2985508}.
Bearing in mind the definition of Pl\"ucker coordinates (see  \eqref{eqMinorsOrd2}, \eqref{eqMinorsOrd3} and \eqref{eqPluckCoord4x4}), it is easy to see
the geometry behind formula \eqref{eqMultDimMAE}: it is nothing but the equation of a hyperplane section of $\LL(n,2n)$, namely the intersection of $\LL(n,2n)$
with the hyperplane
\begin{equation}
 [A:B^{ij}:C^\bullet:\ldots:D^{ij}:D]\in\p(V_{\lambda_n}^*).
\end{equation}
The correct way to formulate the Monge--Amp\`ere test is via the so-called BGG resolution. The same technique provides a similar test for hypersurface
sections of higher degree, that is, with $F$ being a (homogeneous) polynomial  of all the minors of $p_{ij}$ of a certain degree $r>1$. Observe that this
notion of (algebraic) degree has nothing to do with the \emph{order} of the PDE, which is always 2. For instance, {\eqref{eqTheMongeAmpere}} and $p_{11}^2$
are both quadratic in the $p_{ij}$'s, however the former is \emph{linear} in the Pl\"ucker coordinates, whereas the latter is quadratic.  The aforementioned
BGG technique is explained in the paper \cite{MR3603758}.

\subsection{Goursat-type Monge--Amp\`ere equations}\label{secGoursatMA}
A similar expression to  \eqref{eqTheMongeAmpere} describes the so-called  \emph{Goursat-type} (resp., \emph{symplectic Goursat-type}) \emph{Monge--Amp\`ere equation}
\begin{equation}\label{eqDefGoursatMAEs}
 \det(p_{ij}-D_{ij})=0,
\end{equation}
where $D_{ij}$ is a (not necessarily symmetric) $n\times n$ matrix of functions on $M$ (resp., of constants). It is natural to ask oneself whether the class
of (symplectic) Goursat-type Monge--Amp\`ere equations is a proper subclass of the class of (symplectic) Monge--Amp\`ere equations. A straightforward count
of the parameters immediately says yes.
Let us begin with $n=2$. It was already pointed out that the space parametrising the hyperplane sections of $\LL(2,4)$---that is, (symplectic) Monge-Amp\`ere
equations in two variables---is the dual $\p^{4\,\ast}$ of the Pl\"ucker embedding space, see Section \ref{secDualVAr}. On the other hand, the space of
matrices $D_{ij}$ is also 4-dimensional, so that, topological obstruction aside, the two classes of PDEs may well coincide.

Over the complex field, they indeed do.

Over the reals, we have $\E=\{F=0\}$, with
\begin{align}
F=\det(p_{ij}-D_{ij})&=\det(p_{ij})-D_{22}p_{11}+(D_{12}+D_{21})p_{12}-D_{11}p_{22}+\det(D_{ij})\nonumber\\
&=E\det(p_{ij})+Ap_{11}+2Bp_{12}+Cp_{22}+\Delta\label{eqGoursatMA}
\end{align}
and, independently on $D$, the equation $\E$ is always non-elliptic since\footnote{Here we have employed the same notation used in \cite[p.~588]{MR1334205}
for the definition of elliptic/parabolic/hyperbolic Monge--Amp\`ere equations in two dimensions---up to the replacement of symbols $\Delta \leftrightarrow D$.
Observe that the inequality \eqref{eqEqPezzottataDalVitellone} becomes an equality (that is, $\E$ is parabolic) if and only if the matrix $D$ is symmetric.}
\begin{equation}\label{eqEqPezzottataDalVitellone}
(AC-\Delta E-B^2)=-\left(\frac{D_{12}-D_{21}}{2}\right)^2\leq0.
\end{equation}
In other words, for $n=2$, the subclass of  Goursat-type Monge--Amp\`ere equations coincides with the \emph{open} subclass of non-elliptic Monge--Amp\`ere
equations. For $n=3$ a simple dimension count shows that this is no longer possible: the class of (symplectic) Goursat-type Monge--Amp\`ere equations is
($3\cdot 3=9$)-dimensional, whereas all (symplectic) Monge--Amp\`ere equations are parametrised by $\p^{13\ast}$.

\abs-0.1 \aabd-1
From a geometrical standpoint, as we have already stressed in Section \ref{secLagChowTr}, the equation \eqref{eqDefGoursatMAEs} is nothing but the Lagrangian
Chow form of the $(n-1)$-dimensional (linear) variety $\p D$ in $\p(V\oplus V^*)$. Then we are just saying that, in general, not all the hyperplane sections
are the Lagrangian Chow transform of a linearly embedded $\p^{n-1}$ inside the projectivised symplectic space.

The reader should be aware of the fact that,  in the general context of second order PDEs---i.e., when the  hypothesis of being symplectic has been
dropped---the   $n\times n$ matrix $D$ is allowed to depend on the point of $M$. In other words, the   $n$-dimensional subspace $D$ (which we keep denoting
by the same symbol $D$) is actually an $n$-dimensional subdistribution of the contact distribution on $M$. It is $\Cont(M)$-equivariantly associated to the
equation \eqref{eqDefGoursatMAEs} itself. The context will always make it clear, whether $D$ is a distribution or an $n$-dimensional subspace.

Recall that, for $n=2$, the two double fibration pictures \eqref{eqDiagDoubFibrNMenoUno} and \eqref{eqDiagDoubFibrUno} coincide. Then   the Lagrangian Chow
transform can be ``inverted'' simply by taking the characteristic lines (which for $n=2$ are the same as hyperplanes). More precisely, given any
$\E\subset\LL(2,4)$, one defines
\begin{equation}\label{eqPrimaDefIcsE}
X_\E:=\{\Sigma\in\p(V\oplus V^*)\mid \Sigma^{(1)}\textrm{ is tangent to }\E\textrm{ in some point}\}.
\end{equation}
In other words, $X_\E$ is the union of all the characteristics of $\E$ in all its points. More precisely, for any $L\in\E$ consider the characteristic variety
$\Xi_L(\E)$: the points of the latter are, by definition, hyperplanes in $L$, that is lines in $L$. But $L$ is contained into $V\oplus V^*$, so that lines
in $L$ are also lines in $V\oplus V^*$, that is points of $\p(V\oplus V^*)$. So, definition \eqref{eqPrimaDefIcsE} can be rephrased as
\begin{equation}
X_\E=\bigcup_{L\in \E}\Xi_L(\E).
\end{equation}
We stress that the characteristic variety $\Xi(\E)$ is a bundle over $\E$, whereas $X_\E$ is a one-dimensional sub-distribution of $\p(\CC)$, that is
a 2-dimensional \emph{conic sub-distribution} of the contact distribution $\CC$ \cite[Section 7]{Buczyski2019}.

If the equation $\E$  is the Goursat-type Monge--Amp\`ere equation associated, according to \eqref{eqDefGoursatMAEs} to the subdistribution $D\subset\CC$, then
\begin{equation}
 X_\E=\p D\cup \p D^\perp,
\end{equation}
where $\ldots^\perp$ means the symplectic orthogonal.

Observe that the equation \eqref{eqDefGoursatMAEs} above remains unchanged if $D$ is replaced by its orthogonal---the matrix counterpart of taking the
symplectic orthogonal. Then it is not an exaggeration to claim that a Goursat-type Monge--Amp\`ere equation $\E$ is unambiguously determined by its
``inverse Lagrangian Chow form'' $X_\E$. Due to the invariance of the framework, the   sub-distribution $X_\E$ of $\CC$   can by all means replace $\E$ in
the treatment of the equivalence problem. This point of view is at the basis of many works about invariants and classification of Goursat-type Monge--Amp\`ere
equations, see, e.g., \cite{MR2503974,MR2805306,MR1334205,CatalanoFerraioliVinogradov2008,Kushner2009,MR2352610}.

It is worth noticing that the analogous construction of $X_\E$ for multidimensional PDEs is slightly more complicated \cite{MR2985508}. The class of PDEs
``that are the Lagrangian Chow transform of their own inverse  Lagrangian Chow transform''---\emph{re\-con\-struct\-able}, for short---contains in fact more
than the Goursat-type Monge--Amp\`ere equations, but it has not yet been explored completely.

\abs0.1 \aabd1
\subsection{Low-dimensional examples}\label{secLowDimExamples}
Let us recall that, for  $n=2$, the  space $\p^{4\ast}$ naturally parametrises hyperplane sections of $\LL(2,4)$, that is (symplectic) two-dimensional
Monge--Amp\`ere equations, see Section \ref{secDualVAr} and Section \ref{secHypMAEs}. Inside  $\p^{4\ast}$ there sits the three-dimensional dual variety
$\LL(2,4)^*$, viz.,
\begin{equation}
\underbrace{\LL(2,4)^*}_{\textrm{parabolic Monge--Amp\`ere}} \subset \underbrace{\p^{4\ast}}_{\textrm{all Monge--Amp\`ere}},
\end{equation}
and it corresponds precisely to the sub-class of parabolic (Goursat-type symplectic) Monge--Amp\`ere equations, see Section \ref{secGoursatMA}. Indeed,
when $D$ is Lagrangian, i.e., symmetric, the symbol of $F$ has a double root, see Section \ref{secNonDeg} and \eqref{eqGoursatMA}. Since we are working over
the reals, between the subset and the whole space there is also the open domain made of non-elliptic Monge--Amp\`ere equations, that is all Goursat-type
Monge--Amp\`ere equations \eqref{eqDefGoursatMAEs}.

For $n=3$ the stratification becomes more interesting, since the dual variety is singular. We are now in position of interpreting \eqref{eqL36isoDualSinglL36}
in terms of Monge--Amp\`ere equations:
\begin{equation}\label{eqStratL36}
\underbrace{\Sing(\LL(3,6)^*)}_{\textrm{parabolic Monge--Amp\`ere}}\subset \underbrace{\LL(3,6)^*}_{\textrm{linearisable Monge--Amp\`ere}}\subset
\underbrace{\p^{13\ast}}_{\textrm{all Monge--Amp\`ere}}.
\end{equation}
The (9-dimensional) domain of Goursat-type Monge--Amp\`ere equations is between the first two strata. A proof of the fact that the 12-dimensional variety
$\LL(3,6)^*$ corresponds to linearisable Monge--Amp\`ere equation can be found in  \cite[Section 3.6]{doi:10.1093/imrn/rnp134} or in
\hbox{\cite[Theorem 1.4]{Russo2019}}.

The cases $n=3$ and $n=4$ are important in that  another class of Monge--Amp\`ere equations, which is trivial for $n=2$, begins to show up. This is the class
of \emph{integrable Monge--Amp\`ere equations $($by the method of hydrodynamic reductions$)$}, which will be briefly explained in Section \ref{SecInt} below.
For $n=3$ these coincide with the linearisable ones. From \eqref{eqStratL36} it follows that there exist non-integrable Monge--Amp\`ere equations.  In fact,
these are the general ones, since they form two open orbits, represented by
\begin{equation}
\det(p_{ij})=1,\quad \det(p_{ij})=\tr(p_{ij}),
\end{equation}
see \cite[Equation (13)]{doi:10.1093/imrn/rnp134}. In the same paper it is proved that the space of integrable Monge--Amp\`ere equations has dimension 21.
Since linearisable (that is the same as integrable) Monge--Amp\`ere equations correspond to the 12-dimensional variety  $\LL(3,6)^*$, pure dimensional
considerations show that there is a lot of integrable second order PDEs that are not of Monge--Amp\`ere type.

The picture begins to change  starting from $n=4$. First of all, integrable  Monge--Amp\`ere equations do not coincide with the linearisable ones,
\cite[Theorems 1.4, 1.5, 1.6]{Russo2019}. Second, the size of the space parametrising integrable second order PDEs does not grow, as a function of $n$, as fast
as the size of the sub-variety of $\LL(n,2n)^*$ parametrizing integrable  Monge--Amp\`ere equations. This simple observation led Ferapontov and  Doubrov to
conjecture, in  2010, that from $n\geq 4$ a (symplectic) integrable second order PDE must be necessarily an (integrable) Monge--Amp\`ere one
\cite[Section 1]{DoubrovFerapontov}. Even though the case $n=4$ has been recently solved by  Ferapontov,   Kruglikov and  Novikov \cite{2017arXiv170708070F},
to date the conjecture is still unanswered in general.

\abs0.15 \aabd1
\subsection{Integrability by the method of hydrodynamic reductions}\label{SecInt}
As the reader may have noticed, our survey of the geometry of Lagrangian Grassmannians, their hypersurfaces and second order PDEs has begun to border
with ongoing research activities and open problems. It is then the appropriate moment to end it. We will just mention a few recent research results and
current projects to the benefit of the most curious readers, see Section \ref{secPart3} below.

Before that, we briefly outline the notion of hydrodynamic integrability, in view of the central role played here by the Ferapontov conjecture. The special
classes of Monge--Amp\`ere equations introduced so far---including the integrable ones---can all be interpreted as suitable $\Sp_{2n}$-invariant subsets
in $\p V_{\lambda_n}^*$. Nevertheless, the notion of hydrodynamic integrability was born originally at the antipodes of Algebraic Geometry in response
to a rather tangible problem, which is worth recalling.

The first historically recorded ``hydrodynamic reduction'' of a three-dimensional  quasi-linear system of PDEs dates back to 1860 and it is due
to Riemann. His paper \cite{Riemann1998} provides a mathematical treatment of a problem of nonlinear acoustics proposed by von Helmholtz---the propagation
of planar air waves. The system of PDEs describing the problem expressed the temperature $T$ as a function of the three independent variables
$\rho,p,v$---density, pressure and velocity. Riemann's method consisted in postulating the existence of solutions depending on \emph{two} auxiliary
independent variables $r$ and $s$, and then solving the so-obtained \emph{reduced system}.

In 1996, a similar method was employed by Gibbons and Tsarev in order to obtain a ``chain of hydrodynamic reductions'' \cite{GIBBONS199619} out of
a famous multidimensional system of PDEs introduced by  Benney in the seventies \cite{SAPM:SAPM197756181}. Unlike Riemann's work, the so-obtained chain
of hydrodynamic reductions \emph{does not} lead to actual solutions of the original system of PDEs, but the fact that each reduction is \emph{compatible}
reflects a (still unspecified) property of integrability of the system itself.

The idea that a multidimensional (system of) PDEs may be called ``integrable'' if the corresponding ``hydrodynamic reductions''---obtained from it by
a suitable (though straightforward) generalization of the original Riemann's method---are compatible finally reached its maturity in the early 2000's thanks
to the works of Ferapontov and his collaborators (see \cite{doi:10.1093/imrn/rnp134} and references therein). They observed that the condition of being
\emph{integrable $($in the sense of hydrodynamic reductions$)$}  singles out a nontrivial subclass in the class of second order symplectic PDEs, that is precisely
the one mentioned in Section \ref{secLowDimExamples}. They also obtained, for three-dimensional systems, an \emph{integrability test}, that is a PDE imposed
on the left-hand side of an unknown symplectic second order PDE, which is satisfied if and only if the unknown PDE is hydrodynamically integrable. Unlike
the aforementioned Monge--Amp\`ere test \eqref{eqMongAmpTestBrutal}, which is a consequence of the general construction of the BGG resolution,\footnote{Due
to obvious limitations, the details of the Monge--Amp\`ere test based on the BGG resolution cannot  be reviewed here. In  Section \ref{SecRThDegreeSection}
above we have sketched the idea behind it, but for a full account of it the reader should consult \cite{MR3603758}.} Ferapontov's method was based
on computer-algebra computations and this is why the integrability test is now known only for small values of~$n$ (to date, only 3 and 4). In terms of these
tests, Ferapontov conjecture may be recast as follows:
\begin{equation}\label{eqConjFera}
(\textrm{Monge--Amp\`ere test})+( \textrm{integrability test})\equiv (\textrm{integrability test})\, \quad \forall n\geq 4.
\end{equation}
It then all boils down to formalise \eqref{eqConjFera} in a framework which is rich and general enough to make it possible elaborate an answer. A promising
technique is based on the $\CO_n$-structure associated with a non-degenerate hypersurface  in $\LL(n,2n)$, see Section \ref{secPart3}.

\abs-0.15 \aabd-1
In order to see what really means for a (symplectic) PDE $\E=\{F(p_{ij})=0\}$ to be integrable in the aforementioned hydrodynamical sense, we need to explain
in detail the   notion of a  \emph{$k$-phase solution} of  $\E$. Since $\E$ is symplectic, we can identify $\E$ with its fibre, that is a hypersurface
in $\LL(n,2n)$. Then $f\in C^\infty(\R^n)$ is a solution to $\E$ if its Hessian matrix, understood as a Lagrangian plane parametrised by points of $\R^n$,
takes values into $\E$ (combine \eqref{eqElleEffeUno} and \eqref{eqPDE_Hirota}).  In the streak of the aforementioned Riemann's original work, we understand
a $k$-phase solution as a solution $f$ which depends on the independent variables $(x^1,\dots,x^n)\in \R^n$ through the auxiliary variables
$(R^1,\dots, R^k)\in\R^k$, in such a way that the coordinate vector fields $\frac{\partial}{\partial R^i}$ have \emph{rank one}, see Section \ref{secRO}.
In terms of commutative diagrams,
\begin{equation}
\gathered
\xymatrix{
\R^k\ar[rr]^U &&**[r] \E\subset \LL(n,2n)\\
&\R^n\, .\ar[ul]^{\boldsymbol{R}}\ar[ur]_{\hess(f)}&
}\endgathered
\end{equation}

It is no coincidence that these $R^i$'s are  called \emph{Riemann invariants}.
A \emph{$k$-phase solution} of $\E$ is precisely a solution $f$ of $\E$ making commutative the above diagram, with
\begin{equation}\label{eqRankOneCoordVect}
\rank\left[U_*\left( \frac{\partial}{\partial R^i}\right)\right]=1,\quad\forall i=1,\dots,k.
\end{equation}
Basically, a PDE $\E$ is declared to be integrable if it possesses ``sufficiently many'' $k$-phase solutions, for all $k$ (even if it suffices to check it
just for $k=2,3$). More precisely, one couples the given equation $\E$ with auxiliary equations expressing the existence of the functions $U$ and
$\boldsymbol{R}$ and, most importantly, encoding the rank-one condition \eqref{eqRankOneCoordVect}. Then $\E$ is declared  to be integrable if the
so-obtained system is \emph{compatible}. More details can be found, e.g.,  in \cite{doi:10.1093/imrn/rnp134}.

It is worth observing that the (physically motivated) notion of a $k$-phase solution corresponds to the purely algebro-geometric concept of a
\emph{$k$-secant variety}. This interesting parallel is the main motivation behind the recent work of Russo  \cite{Russo2019}.

\subsection{A selection of recent research results}\label{secPart3}

The main consequence of the non-degeneracy of a hypersurface $\E$ in $\LL(n,2n)$ is the presence of a $\CO_n$-structure on $\E$. This is essentially due to
the reduction of $\GL(V)$, the zero-degree part of $P$, to $\CO(\Smbl(\E))$, the subgroup of $\GL(V)$ preserving the line $\Smbl(\E)$, see definition
\eqref{eqDefSimbProject}. Because the $P$-principal bundle $\Sp_{2n}\longrightarrow \LL(n,2n)$ is made of second-order frames (see Section \ref{sec2NdFrBund}),
such a  $\CO_n$-structure on $\E$ is \emph{not} a conformal metric. This makes things even more intriguing.

There is not yet in the literature a systematic treatment of such $\CO_n$-structures and this is not the appropriate place to start one. Worth to mention
however is the skilful work \cite{MR2765729} by Smith, where $n$ is assumed to be 3 and hence $\CO_3$-structures are the same as $\GL_2$-structures.
There the author  even goes beyond the class of hypersurfaces (5-folds) in $\LL(3,6)$, and studies the equivalence problem of arbitrary $\GL_2$-structures
in dimension~5. Invariants are extracted from a preferred  principal connection which is associated with each such structure. In particular, he finds the
embeddability conditions (i.e., those ensuring that an abstract $\GL_2$-structure in dimension 5 can be realised as a hypersurface in $\LL(3,6)$) and lists
several non-equivalent classes of second order symplectic PDEs in three independent variables.

An analogous treatment of the 4-dimensional case is still lacking in the literature. Nevertheless it is worth to mention the recent preprint by {Ferapontov},
{Kruglikov} and {Novikov}, who answered  the Ferapontov conjecture for $n=4$ \cite{2017arXiv170708070F}.

Concerning the Lagrangian Chow form and the correspondence between substructures of the contact distribution and second order PDEs, it is worth to mention
the work \cite{The2018} by The and the almost simultaneous work \cite{AGMM_CommContMath} by the authors and  Alekseevsky. The problem dealt with there is
that of constructing a PDE admitting a prescribed simple (complex) Lie group of symmetries. The departing point is the so-called \emph{sub-adjoint variety}
\cite[Section 8]{Buczyski2019} of a rational homogeneous contact manifold, which is an example of a \emph{conic sub-distribution} of the contact distribution,
see Section \ref{secLagChowTr}.  Besides these highly symmetric cases, there is still no systematic treatment of ``higher degree'' analogues of Goursat-type
Monge--Amp\`ere equations.

Especially to the reader who is wondering ``why always second order'' we may suggest the paper \cite{MannoMoreno2016} where an analogous approach to the one
proposed here has been applied to the natural third order analogues of Monge--Amp\`ere equations.

\section{Appendix: a guide to reading this volume}\label{secPart3.1}

The present paper was entirely dedicated to the geometry of the Lagrangian Grassmannian and its hypersurfaces. However, it should not be forgotten
that the Lagrangian Grassmannian bundle $M^{(1)}$ over a contact manifold $M$ is but an example of the \emph{variety of integral elements} of an Exterior
Differential System (EDS). The theory of EDS'es represents one of the most general  frameworks for studying (system of) PDEs from the point of view of
differential geometry (an alternative approach is based on jet spaces \cite{MR861121}). It was born with the pioneering works of Pfaff
\cite{pfaff1818methodus},  Frobenius and  Darboux \cite{Darboux1882} and Cartan \cite{Cartan_EDS}. Later it was perfected through many contributions.
All details about the modern incarnation of theory can be found in the excellent book \cite{MR1083148} by Bryant, Chern, Gardner, Goldschmidt and Griffiths.
However, the size of the volume may be discouraging for those who seek a swift and workable introduction to the topic.  McKay's paper \cite{McKay_BCP} serves
precisely such a purpose.

Smith's paper \cite{Smith2019} is, in a sense, complementary to McKay's one. While the latter is concerned with differential ideals and PDEs,  the former
focuses instead on the geometry of the set of integral elements of an EDS, understood as a sub-bundle of the Grassmannian bundle. The vertical bundle of these
sub-bundles, known as \emph{tableaux}, the characteristic variety and its incidence correspondence are all examined in detail.

\abs-0.1
The paper \cite{Russo2019} by Russo is a natural companion to the works by Ferapontov and his collaborators on the geometry of hydrodynamic integrability
\cite{DoubrovFerapontov,doi:10.1093/imrn/rnp134}. The author carefully explains several algebro-geometric notions and theorems that are made use of,
more or less explicitly,  in Ferapontov's works. In particular, some classical results on the geometry of secant varieties are reviewed, and a nice technique
is employed to deal with the cases of $\LL(3,6)$ and $\LL(4,8)$: the analogy of these cases with the twisted cubic in~$\p^3$ and the rational normal curve
in $\p^4$, respectively.

The notion of a contact structure does not pertain exclusively to the realm of differentiable manifolds. The   parallel idea of a \emph{complex contact
manifold} is reviewed in \cite{Buczyski2019}, by Buczy\'nski and one of us (Moreno). However,  the main purpose of the paper is that of underlying important
and unexpected bridges between the complex-analytic and the real-differentiable setting. In particular, there are discussed the twistor correspondence for
quaternion-K\"ahler manifolds and certain substructures of the contact distribution that can be studied via the Cartan's method of equivalence.

Panasyuk's paper \cite{Panasyuk2019} reviews an interesting correspondence, basically due to  Gelfand and Zakharevich, between the notion of a
bi-Hamiltonian system and the notion of a \emph{Veronese web}. The former is a powerful tool, widely exploited in the theory of integrable systems, that are
PDEs admitting a particularly rich and well-behaved set of (higher) symmetries and/or conservation laws. The latter is a purely geometric construction,
generalising that of a web: it is a family, parametrised by  $\p^1$, of foliations such that the annihilators describe a rational normal curve in the
projectivised cotangent bundle. Such a geometric interpretation adds some clarity to the integrable systems' area of research, which features many excellent
techniques but sometimes lacks theoretical rigor.

The reader who was surprised by the identification of $\LL(2,4)$ with the Lie quadric may appreciate  Jensen's short review of Lie sphere geometry
\cite{Jensen2019}, which deals with the space $\mathcal{S}(\R^3)$ of generalised spheres  in the Euclidean space $\R^3$. ``Generalised'' means that it
encompasses the   spheres of zero radius (points) and those of infinite radius (planes)   as well. Consider the pseudo-Euclidean vector space
$V = \R^{4,2}$    of signature $(4,2)$.  Then $\mathcal{S}(\R^3)$ is identified with the points of the Lie quadric $Q = \p V_0$, that is  the
projectivisation of the isotropic cone $V_0$ in $V$. The author studies   the quadric as a homogeneous space $Q = \OOO(4, 2)/P$ where $P$ is the parabolic
subgroup that stabilises an isotropic line. Lines in $Q$ correspond to two-dimensional absolutely isotropic $2$-planes in $V$.

Finally,  Musso and Nicolodi's paper \cite{Musso2019} provides a lucid introduction to Laguerre geometry with a clean presentation of the
fundamental constructions. It contains helpful comparisons to surface theory in other, classical geometries. Its subject represents a perfect arena to show
the potential of the standard method of moving frames and of EDS'es.  The reader will be pleased to see how geometry  adds  some perspective and helps to
demystify the more technical aspects of the Cartan--K\"ahler theorem as well as of the frame adaptation.

\bibliographystyle{plainurl}
\bibliography{BibUniver}

\end{document}